\documentclass[12pt,a4paper]{article}
\usepackage{mathrsfs}
\usepackage{amssymb}
\usepackage{amsmath}
\usepackage{amsfonts}
\usepackage{amstext}
\usepackage{amsthm}
\usepackage{graphicx}
\usepackage[top=1in, bottom=1in,left=0.85in,right=0.85in]{geometry}

\def \qed {\hfill \vrule height6pt width 6pt depth 0pt}

\begin{document}

\title{Two Whyburn type topological theorems and its applications to Monge-Amp\`{e}re equations
\thanks{Research supported by NNSF of China (No. 11261052, No. 11201378).}}
\author{{\small  Guowei Dai\thanks{Corresponding author. \newline
\text{\quad\,\, E-mail address}: daiguowei@nwnu.edu.cn.
}
} \\
{\small Department of Mathematics, Northwest Normal
University, Lanzhou, 730070, P.R. China}\\
}
\date{}
\maketitle

\begin{abstract}
In this paper we correct a gap of Whyburn type topological lemma and establish two superior limit theorems.
As the applications of our Whyburn type topological theorems, we study the following Monge-Amp\`{e}re equation
\begin{eqnarray}
\left\{
\begin{array}{lll}
\det\left(D^2u\right)=\lambda^N a(x)f(-u)\,\, &\text{in}\,\, \Omega,\\
u=0~~~~~~~~~~~~~~~~~~~~~~\,\,&\text{on}\,\, \partial \Omega.
\end{array}
\right.\nonumber
\end{eqnarray}
We establish global bifurcation results for the problem.
We find intervals of $\lambda$ for the existence, multiplicity and nonexistence of strictly convex solutions for this problem.
\\ \\
\textbf{Keywords}: Superior limit; Topological methods; Monge-Amp\`{e}re equations; Bifurcation; Strictly convex solutions
\\ \\
\textbf{MSC(2000)}: 12J99; 34C23; 34D23; 34L05; 35J60
\end{abstract}\textbf{\ }

\numberwithin{equation}{section}

\numberwithin{equation}{section}

\section{Introduction}

\quad\, The Whyburn's limit theorem [\ref{Whyburn}, Theorem 9.1] is an important tool in the study of differential equations theory, see for example, [\ref{ACD}, \ref{AT}, \ref{ADT}, \ref{Dai}] and references cited therein. However, if the collection of the infinite sequence of sets
is unbounded, the Whyburn's limit theorem cannot be used directly because the collection may not be relatively compact.
In order to overcome this difficulty, the authors of [\ref{LiSun}] and [\ref{MaAn}] (independently) establish
the following topological lemma:
\\ \\
\textbf{Lemma 1.1.} \emph{Let $X$ be a Banach space and let $\left\{C_n\right\}$ be a family of closed connected subsets of $X$. Assume that:}
\\

(a) \emph{there exist $z_n\in C_n$, $n=1,2,\ldots$, and $z^*\in X$, such that $z_n\rightarrow z^*$;}

(b) \emph{$r_n=\sup \left\{\Vert x\Vert: x\in C_n\right\}=+\infty$;}

(c) \emph{for every $R>0$, $\left(\cup_{n=1}^{+\infty} C_n\right)\cap B_R$ is a relatively compact set of $X$, where}
\begin{equation}
B_R=\{x\in X:\Vert x\Vert\leq R\}.\nonumber
\end{equation}
\noindent \emph{Then there exists an unbounded component $C$ in $D =\limsup_{n\rightarrow +\infty}C_n$ and $z^*\in \mathfrak{C}$.}
\\

Lemma 1.1 has been used several times in
the literature to analyze the global structure of nontrivial solutions in wide classes of elliptic boundary value problems for equations and
systems, see for example, [\ref{DML}, \ref{DMW}, \ref{Jankowski},  \ref{MGX}, \ref{XOC}, \ref{XSO}]. Unfortunately, the proof of Lemma 1.1 contains a gap.
In order to prove the component $C$ containing $z^*$ is unbounded, they assume $C$ is bounded then get a contradiction.
Indeed, if $\left\Vert z^*\right\Vert=+\infty$, $C$ must be unbounded. So the assumption of $C$ being bounded is impossibly holding in this case.
The first aim of this paper is to correct this gap and establish two superior limit theorems, which will be done in Section 2.

As applications of our superior limit theorems, we study convex solutions of the Monge-Amp\`{e}re equations.
The Monge-Amp\`{e}re equations are a type of important fully nonlinear elliptic equations [\ref{GT}, \ref{Guan0}, \ref{Guan4}, \ref{P3}, \ref{T}]. The study
of the Monge-Amp\`{e}re equations has been received considerable attentions in history, which is motivated by the Minkowski problem [\ref{CY1}, \ref{Minkowski}, \ref{P0}, \ref{P2}] and the Weyl problem [\ref{Guan5}, \ref{Nirenberg}] in classical Euclidean geometry and conformal
geometry. Existence and regularity results of the Monge-Amp\`{e}re equations can be found in [\ref{CNS}, \ref{CY1}, \ref{CY}, \ref{CY2}, \ref{GT}, \ref{Guan3}, \ref{Guan1}, \ref{Guan2},  \ref{K}, \ref{L}, \ref{Nirenberg}, \ref{P}, \ref{P1}, \ref{P2}, \ref{ZW}] and the reference therein.
Consider the following real Monge-Amp\`{e}re equations
\begin{eqnarray}\label{MB}
\left\{
\begin{array}{lll}
\det\left(D^2u\right)=\lambda^N a(x)f(-u)\,\, & \text{in}\,\, B,\\
u=0~~~~~~~~~~~~~~~~~~~~~\,\,& \text{on}\,\, \partial B,
\end{array}
\right.
\end{eqnarray}
where $D^2u=\left(\frac{\partial^2 u}{\partial x_i\partial x_j}\right)$ is the Hessian
matrix of $u$, $B$ is the unit ball of $\mathbb{R}^N$, $a\in C\left(\overline{B}, [0,+\infty)\right)$ is a radially symmetric weighted function and $a(r):=a(\vert x\vert)\not\equiv 0$ on any subinterval of $[0, 1]$, $\lambda$ is a
positive parameter and $f:[0,+\infty)\rightarrow [0,+\infty)$
is a continuous function such that it does not vanish on any entire interval. The study of problem (\ref{MB}) in
general domains of $\mathbb{R}^N$ may be found in [\ref{CNS}, \ref{GT}]. Kutev
[\ref{K1}] investigated the existence of strictly convex radial solutions
of problem (\ref{MB}) when $f(s)=s^p$ and $a(x)\equiv 1$. Delano [\ref{Del}] treated the
existence of convex radial solutions of problem (\ref{MB}) for a class of
more general functions, namely $\lambda\exp f(\vert x\vert,u,\vert \nabla u\vert)$.
Guan and Lin [\ref{Guan6}] studied more general Monge-Amp\`{e}re equations with weighted nonlinearity.
As in [\ref{HW}, \ref{K1}], we can easily show that the radially symmetric solutions of problem (\ref{MB})
satisfies the following boundary value problem
\begin{equation}\label{MO}
\left\{
\begin{array}{l}
\left(\left(u'\right)^N\right)'=\lambda^NNr^{N-1} a(r)f(-u),\,\,  r\in(0,1),\\
u'(0)=u(1)=0.
\end{array}
\right.
\end{equation}
By a solution of problem (\ref{MO}) we understand that it is a
function which belongs to $C^2[0,1]$ and satisfies (\ref{MO}). It is easy to verify that any negative solution
of problem (\ref{MO}) is strictly convex in $(0,1)$. Wang [\ref{W}], Hu-Wang [\ref{HW}] established several
criteria for the existence, multiplicity and nonexistence of strictly convex solutions for problem (\ref{MO}) with $a(r)\equiv 1$ by using
fixed index theorem. However, there is no any information on the bifurcation points and the
optimal intervals for the parameter $\lambda$ so as to ensure existence of single or multiple convex solutions.

We shall use nonlinear analysis method and our superior limit theorems to find the optical intervals for
the parameter $\lambda$ so as to ensure existence of single or multiple solutions and possible bifurcation points; thus extend and improve
the corresponding results of [\ref{HW}, \ref{W}]. To do this, in Section 4, we shall establish a global bifurcation theorem
for problem (\ref{MO}) with $f(u)=u^N+g(u)$, i.e.,
\begin{equation}\label{Mg}
\left\{
\begin{array}{l}
\left(\left(u'\right)^N\right)'=\lambda^NNr^{N-1} a(r)\left(\left(-u\right)^N+g(-u)\right),\,\, r\in(0,1),\\
u'(0)=u(1)=0,
\end{array}
\right.
\end{equation}
where $g:[0,+\infty)\rightarrow [0,+\infty)$ satisfies $\lim_{ s\rightarrow0^+}g(s)/s^N=0$.
Concretely, we shall show that $\left(\lambda_1,0\right)$ is a bifurcation point of problem (\ref{Mg}) and there exists an
unbounded continuum of convex solutions, where $\lambda_1$ is the first eigenvalue of problem (\ref{Mg})
with $g\equiv 0$. Note that problem (\ref{Mg}) is a type of nonlinear equation. Hence, the common index formula
involving linear map cannot be used here. In order to overcome this difficulty, we shall study an auxiliary
eigenvalue problem in Section 3, which has an independent interest, and establish an
index formula for it. Then by use of the index formula of the auxiliary problem, we prove an index formula
involving problem (\ref{Mg}) which guarantees $\left(\lambda_1,0\right)$ is a bifurcation point of nontrivial solutions
of problem (\ref{Mg}).

On the basis of the above global bifurcation results and superior limit theorems, we investigate the existence of strictly convex solutions of problem (\ref{MO}). We shall give the optimal intervals for the parameter $\lambda$
so as to ensure existence of single or multiple strictly convex solutions in Section 5 and 6. Our results extend and/or improve
the corresponding results of [\ref{HW}, \ref{L1}, \ref{W}].
Basis on the results on unit ball, we also study problem (\ref{MB}) on a general domain $\Omega$ in Section 7, i.e.,
\begin{eqnarray}\label{MAh1}
\left\{
\begin{array}{lll}
\det\left(D^2u\right)=\lambda^N a(x)f(-u)\,\, &\text{in}\,\, \Omega,\\
u=0~~~~~~~~~~~~~~~~~~~~~\,\,&\text{on}\,\, \partial \Omega,
\end{array}
\right.
\end{eqnarray}
where $\Omega$ is a bounded convex domain of $\mathbb{R}^N$ with smooth boundary and $0\in\text{Int}\Omega$.
It is well known [\ref{GT}] that problem (\ref{MAh1}) is elliptic only when the Hessian matrix $D^2u$ is positive definite and
it is therefore natural to confine our attention to convex solutions, nonnegative functions $f$ with $f(s)>0$ for $s>0$ and nonnegative weighted function $a(x)$ with $a(x)>0$ for $x\in \Omega$.
Obviously, any convex solution of problem (\ref{MAh1}) is negative and strictly convex. In [\ref{ZW}], the authors have proved a lemma concerning the comparison between domains for problem (\ref{MAh1}) with $f(s)=e^s$ and $a(x)\equiv 1$ by sub-supersolution method. We shall show that this lemma is also valid for problem (\ref{MAh1}).
Using this domain comparison lemma and the results on unit ball, we can prove some existence and nonexistence of
convex solutions for problem (\ref{MAh1}).

\section{Superior limit theorems}

\quad\, In this section, we amend the above mentioned gap by establishing the following two superior limit theorems.
\\ \\
\textbf{Theorem 3.1.} \emph{Let $X$ be a normal space and let $\left\{C_n\right\}$ be a sequence of unbounded connected subsets of $X$. Assume that:}
\\

(i) \emph{there exists $z^*\in \liminf_{n\rightarrow +\infty} C_n$ with $\left\Vert z^*\right\Vert<+\infty$;}

(ii) \emph{for every $R>0$, $\left(\cup_{n=1}^{+\infty} C_n\right)\cap B_R$ is a relatively compact set of $X$, where}
\begin{equation}
B_R=\left\{x\in X:\left\Vert x\right\Vert\leq R\right\}.\nonumber
\end{equation}
\noindent \emph{Then $D=\limsup_{n\rightarrow +\infty}C_n$ is unbounded closed connected.}
\\ \\
\textbf{Proof.} Let $X_R=X\cap B_R$ for any $R>0$. Then $X_R$ is a metric subspace under the induced topology of $X$. Let $A_n=C_n\cap B_R$. Clearly, we have
$\cup_{n=1}^{+\infty} A_n=\left(\cup_{n=1}^{+\infty} C_n\right)\cap B_R$.
So $\cup_{n=1}^{+\infty} A_n$ is relatively compact in $X_R$. Furthermore, $z^*\in \liminf_{n\rightarrow +\infty} C_n$ implies that every
neighborhood $U\left(z^*\right)$ of $z^*$ contains points of all but a finite number of the sets
of $\left\{C_n\right\}$. So there exists a positive integer $N$ such that for $n >N$, $U\left(z^*\right)\cap C_n\neq\emptyset$.
Since $\left\Vert z^*\right\Vert<+\infty$, we can take $R>0$ large enough such that $U\left(z^*\right)\subseteq B_R$. Thus, $U\left(z^*\right)\cap A_n=U\left(z^*\right)\cap C_n\neq\emptyset$ for $n >N$.
So we have $z^*\in \liminf_{n\rightarrow +\infty} A_n$.
By Theorem 9.1 of [\ref{Whyburn}], it follows that $A=\limsup_{n\rightarrow +\infty}A_n$ is connected in $X_R$.

We claim that $B:=\left(\limsup_{n\rightarrow +\infty}C_n\right)\cap B_R=A$.
For $x\in A$, then any neighborhood $V$ in $X_R$ of $x$ contains contains points of
infinitely many sets of $\left\{C_n\cap B_R\right\}$. So there exist $x_{n_i}\in C_{n_i}\cap B_R$ such that
$x_{n_i}\rightarrow x$ as $i\rightarrow+\infty$. It follows that $x\in B_R$ and $x\in \limsup_{n\rightarrow +\infty}C_n$, i.e., $x\in B$.
Conversely, if $x\in B$, any neighborhood $V$ in $X_R$ of $x$ contains a point $z$ of $\limsup_{n\rightarrow +\infty}C_n$ and thus $V$, a neighborhood of $z$, contains points of infinitely many of the sets of $C_n\cap B_R$.
It follows that $x\in A$. Hence, $B$ is connected. By the arbitrary of $R$, we get that $D$ is connected.
From [\ref{Whyburn}], we know that $D$ is closed.

Next, we show that $D$ is unbounded.
Suppose on the contrary that $D$ is bounded. It is easy to see that $D$ is a compact set of $X$ by (ii) and the fact of $z^*\in D$.
Let $U_\delta$ be a $\delta$-neighborhood of $D$. So we have that
\begin{equation}\label{Why}
\partial U_\delta\cap D=\emptyset.
\end{equation}
By (i) and the connectedness of $C_n$, there exists an integer $N_0 >0$, such that for all
$n>N_0$, $C_n\cap \partial U_\delta\neq \emptyset$. Take $y_n \in C_n\cap \partial U_\delta$, then $\left\{y_n:n>N_0\right\}$ is
a relatively compact subset of $X$, so there exist $y^*\in \partial U_\delta$ and a subsequence $\left\{y_{n_k}\right\}$ such
that $y_{n_k}\rightarrow y^*$. The definition of superior limit shows that $y^*\in D$. Therefore, $y^*\in \partial U_\delta\cap D$. However, this
contradicts (\ref{Why}). Therefore, $D$ is unbounded.\qed
\\ \\
\textbf{Theorem 2.2.} \emph{Let $X$ be a normal space and let $\left\{C_n\right\}$ be a sequence of subsets of $X$. Assume that
there exists $z^*\in \liminf_{n\rightarrow +\infty} C_n$ with $\left\Vert z^*\right\Vert=+\infty$.
Then there exists an unbounded component $C$ in $E:=\liminf_{n\rightarrow +\infty} C_n$ and $z^*\in C$.}
\\ \\
\textbf{Proof.} Let $C$ denote the component of $E$ containing $z^*$. The unboundedness of $z^*$ shows that $C$ is unbounded.\qed\\

Clearly, the assumptions of Theorem 2.2 are weaker than the corresponding ones of Lemma 1.1 in the case of $\left\Vert z^*\right\Vert=+\infty$; the conclusion of Theorem 2.1 is stronger than that of Lemma 1.1 in the case of $\left\Vert z^*\right\Vert<+\infty$.
So we not only correct the above mentioned gap but also improve the results of Lemma l.1. Moreover, the following example shows that $D$ and $E$ may not be connected under the assumptions of Theorem 2.2.
\\
\\
\textbf{Example 2.1.} Let $X=\mathbb{R}\cup\{+\infty,-\infty\}$. Consider sequence $\left\{A_n\right\}$ in $X$ as the following
\begin{equation}
A_1=[2,+\infty]\cup [2,+\infty], \nonumber
\end{equation}
\begin{equation}
A_{2m}=\left[0,1+\frac{1}{2m}\right]\cup [3,+\infty],\,\, m=1,2,3,\cdots,\nonumber
\end{equation}
\begin{equation}
A_{2m+1}=\left[0,2-\frac{1}{2m+1}\right]\cup [3,+\infty],\,\, m=1,2,3,\cdots.\nonumber
\end{equation}
Let $z^*=+\infty$. Clearly, $z^*\in \liminf_{n\rightarrow +\infty} A_n$ with $\left\Vert z^*\right\Vert=+\infty$.
While, it is not difficult to show that
\begin{equation}
\limsup_{n\rightarrow +\infty} A_{n}=\left[0,2\right]\cup [3,+\infty],\,\,\liminf_{n\rightarrow +\infty} A_{n}=\left[0,1\right]\cup [3,+\infty].\nonumber
\end{equation}
Clearly, $\limsup_{n\rightarrow +\infty} A_{n}$ and $\liminf_{n\rightarrow +\infty} A_{n}$ are not connected. \qed

\section{An auxiliary eigenvalue problem}

\quad\, Consider the following auxiliary problem
\begin{equation}\label{AP}
\left\{
\begin{array}{l}
-\left(\left\vert v'(r)\right\vert^{p-2}v'(r)\right)'=
\mu^{p-1} (p-1)r^{p-2}a(r)\vert v(r)\vert^{p-2}v(r),\,\, r\in(0,1),\\
v'(0)=v(1)=0,
\end{array}
\right.
\end{equation}
where $p\in[2,+\infty)$.
Let $X$ be the Banach space $C[0,1]$ with the norm
\begin{equation}
\Vert v\Vert=\sup_{r\in[0,1]}\vert v(r)\vert.\nonumber
\end{equation}
Define the map $T_\mu^p:X\rightarrow X$ by
\begin{equation}
T_\mu^pv=\int_1^r\varphi_{p'}\left(\int_s^0 \mu^{p-1}
(p-1)\tau^{p-2}a(\tau)\varphi_p(v)\,d\tau\right)\,ds,\,\, 0\leq r\leq1,\nonumber
\end{equation}
where $\varphi_p(s)=\vert s\vert^{p-2}s$, $p'=p/(p-1)$.
It is easy to see that $T_\mu^p$ is continuous and compact.
Clearly, problem (\ref{AP}) can be equivalently written as
\begin{equation}
v=T_\mu^p v.\nonumber
\end{equation}
\indent Firstly, we have the following existence and uniqueness result for problem (\ref{AP}).
\\ \\
\textbf{Lemma 3.1.} \emph{If $(\mu, v)$ is a solution of (\ref{AP})
and $v$ has a double zero, then $v \equiv 0$.}
\\ \\
\textbf{Proof.} Let $v$ be a solution of problem (\ref{AP}) and $r_*\in[0, 1]$ be a double zero.
We note that $v$ satisfies
\begin{equation}
v(r)=\int_{r_*}^r\varphi_{p'}\left(\int_s^{r_*} (p-1)\mu^{p-1}\tau^{p-2}a(\tau)\varphi_p(v)\,d\tau\right)\,ds.\nonumber
\end{equation}
Firstly, we consider $r\in\left[0, r_*\right]$. Then we have
\begin{eqnarray}
\vert v(r)\vert&\leq&\varphi_{p'}\left(\int_{r}^{r_*}(p-1)
\mu^{p-1}\tau^{p-2}a(\tau)\varphi_p(\vert v\vert)\,d\tau\right).\nonumber
\end{eqnarray}
It follows from above that
\begin{eqnarray}
\varphi_p(\vert v\vert)&\leq&\mu^{p-1}\int_{r}^{r_*}(p-1)\tau^{p-2}a(\tau)\varphi_p(\vert v\vert)\,d\tau. \nonumber
\end{eqnarray}
By the Gronwall's inequality, we get $v \equiv 0$ on $\left[0, r^*\right]$.
Similarly, we can get $v \equiv 0$ on $\left[r^*, 1\right]$
and the proof is completed.\qed\\

Set $W_c^{1,p}(0,1):=\left\{v\in W^{1,p}(0,1): v'(0)=v(1)=0\right\}$ with the norm
\begin{equation}
\Vert v\Vert_w=\left(\int_0^1\vert v'\vert^p\,dr\right)^{1/p}.\nonumber
\end{equation}
Then it is easy to verify that $\Vert \cdot\Vert_w$ is the equivalent norm of $W^{1,p}(0,1)$; hence $W_c^{1,p}(0,1)$ is a real Banach space.
\\ \\
\textbf{Definition 3.1.} We call that $v\in W_c^{1,p}(0,1)$ is the weak solution of
problem (\ref{AP}), if
\begin{equation}
\int_{0}^1\left\vert v'\right\vert^{p-2}v'\phi'\,{d}r=
(p-1)\mu^{p-1}\int_{0}^1r^{p-2}a(r)\vert v\vert^{p-2}v\phi\,{d}r\nonumber
\end{equation}
for any $\phi\in W_c^{1,p}(0,1)$.\\
\\
\indent For the regularity of weak solution, we have the following result.\\
\\
\textbf{Lemma 3.2.} \emph{Let $v$ be a weak solution of problem (\ref{AP}), then $v$ satisfies problem (\ref{AP}).}\\
\\
\textbf{Proof.} According to Definition 3.1, we have
\begin{equation}
-\left(\left\vert v'(r)\right\vert^{p-2}v'(r)\right)'=
\mu^{p-1}(p-1)r^{p-2}a(r)\vert v(r)\vert^{p-2}v(r)\,\,\text{in}\,\,(0,1)\nonumber
\end{equation}
in the sense of distribution, i.e.,
\begin{equation}
-\left(\left\vert v'(r)\right\vert^{p-2}v'(r)\right)'
=\mu^{p-1}(p-1)r^{p-2}a(r)\vert v(r)\vert^{p-2}v(r)\,\,\text{in}\,\,(0,1)\setminus I\nonumber
\end{equation}
for some $I\subset (0,1)$ which satisfies $\text{meas}\{I\}=0$.
Furthermore, by virtue of the compact embedding of $W_c^{1,p}(0,1)\hookrightarrow
C^{\alpha}[0,1]$ with some $\alpha\in(0,1)$ (see [\ref{E}]), we obtain that $v\in C^{\alpha}[0,1]$.
Thus, we have that $\lim_{r\rightarrow r_0} \mu^{p-1}(p-1)r^{p-2}a(r)\vert v(r)\vert^{p-2}v(r)$ exists for any $r_0\in I$.
Letting $u:=-\varphi_p\left(v'\right)$, we have
\begin{equation}
\lim_{r\rightarrow r_0} u'(r)=\lim_{r\rightarrow r_0} \mu^{p-1}(p-1)r^{p-2}a(r)\vert v(r)\vert^{p-2}v(r).\nonumber
\end{equation}
The above relation follows that
$\lim_{r\rightarrow r_0} u'(r)$ exists for any $r_0\in I$. Thus,
Proposition 1 of [\ref{DMW}] follows that $u\in C^1(0,1)$,
which implies that $v$ satisfies problem (\ref{AP}).\qed\\

\indent Define the functional $J$ on $W_c^{1,p}(0,1)$ by
\begin{equation}
J(v)=\int_0^1\frac{1}{p}\left\vert v'(r)\right\vert^p\,dr
-\mu^{p-1}\frac{p-1}{p} \int_0^1r^{p-2}a(r)\vert v\vert^p\,dr.\nonumber
\end{equation}
It is not difficult to verify that the critical points of $J$ are the weak solutions of problem (\ref{AP}).
Taking $f_1(v):=\int_0^1\frac{1}{p}\left\vert v'(r)\right\vert^p\,dr$ and
$f_2(v):=\frac{p-1}{p}\int_0^1r^{p-2}a(r)\vert v\vert^p\,dr$,
consider the following eigenvalue problem
\begin{equation}\label{APE}
A(v)=\eta B(v),
\end{equation}
where $A=\partial f_1$ and $B=\partial f_2$ denote the sub-differential of $f_1$ and $f_2$, respectively
(refer to [\ref{C}] for the details of sub-differential).

By some simple computations, we can show that there exists positive $\delta$ such that
\begin{equation}\label{r1}
\frac{f_1(v)}{f_2(v)}\geq \delta
\end{equation}
for any $v\in W_c^{1,p}(0,1)$ and $v\not \equiv0$.
Moreover, we have the following result.
\\ \\
\textbf{Lemma 3.3.} \emph{Put $\eta_1(p)=\inf_{v\in W_c^{1,p}(0,1),v\not\equiv 0}\frac{f_1(v)}{f_2(v)}$.
Then we have that} \\

\emph{(i) (\ref{APE}) has no nontrivial solution for $\eta\in\left(0,\eta_1(p)\right)$;}

\emph{(ii) $\eta_1(p)$ is simple, i.e., (\ref{APE}) has a positive solution and the set of all solutions of (\ref{APE}) is an one dimensional linear subspace of $W_c^{1,p}(0,1)$;}

\emph{(iii) (\ref{APE}) has a positive solution if and only if $\eta=\eta_1(p)$.}\\
\\
\textbf{Proof.} Let $W_c^{1,p}(0,1)=:V$. We denote by $\Phi(V)$ the family of all proper lower semi-continuous convex functions $\varphi$
from $V$ into $(-\infty,+\infty]$, where ``proper'' means that the effective domain $D(\varphi)=\{x\in V: \varphi(x)<+\infty\}$ of $\varphi$ is not empty.

Next, we verify the conditions (A0)--(A4) of [\ref{IO}]. Clearly, we have that $f_1$, $f_2\in \Phi(V)$, $D(f_1)=D(f_2)=V$ and $V\subset L_{\text{loc}}^1(0,1)$,
i.e., condition (Al) is satisfied (by taking $\Omega=(0,1)$). Let $R(v):=f_2(v)/f_1(v)$. Then we have $R(\vert v\vert)\geq R(v)$ for $\forall v\in V$.
It is easy to see that $f_1(v)\geq 0$ for $\forall v\in V$ and $f_1(v)=0$ if and only if $v=0$. Note that (\ref{r1}) implies
$\exists u\in V$ s.t. $u\neq 0$ and $R(u)=\sup\{R(v); v\in V,v\neq 0\}$. So condition (A2) is verified. Taking $\alpha=p$,
we have $f_i(tv)=t^\alpha f_i(v)$ for $\forall v\in V^+=\{w\in V; w(r)\geq0\,\, \text{a.e.}\,\, r\in(0,1)\}$, $\forall t>0$, $i=1,2$.
Thus, condition (A3) is satisfied. For any $u$, $v\in V^+$, we define $(u\vee w)(r)=\max(u(r),w(r))$, $(u\wedge w)(r)=\min(u(r), w(r))$,
$I_1=\{r\in[0,1]:u(r)\geq w(r)\}$ and $I_2=\{r\in[0,1]:u(r)<w(r)\}$. Then we have
\begin{eqnarray}
f_1(u\vee w)+f_1(u\wedge w)&=&\int_0^1\frac{1}{p}\left\vert (u\vee w)'(r)\right\vert^p\,dr+\int_0^1\frac{1}{p}\left\vert (u\wedge w)'(r)\right\vert^p\,dr\nonumber\\
&=&\int_{I_1}\frac{1}{p}\left\vert u'\right\vert^p\,dr+\int_{I_2}\frac{1}{p}\left\vert w'\right\vert^p\,dr
+\int_{I_1}\frac{1}{p}\left\vert  w'\right\vert^p\,dr+\int_{I_2}\frac{1}{p}\left\vert  u'\right\vert^p\,dr\nonumber\\
&=&\int_{0}^1\frac{1}{p}\left\vert u'\right\vert^p\,dr+\int_{0}^1\frac{1}{p}\left\vert  w'\right\vert^p\,dr\nonumber\\
&=&f_1(u)+f_1(w)\nonumber.
\end{eqnarray}
Similarly, we can also show that $f_2(u\vee w)+f_2(u\wedge w)=f_2(u)+f_2(w)$. Hence, condition (A4) is verified. Finally, Lemma 3.1 and 3.2 imply that
every nonnegative nontrivial solution $u$ of (\ref{APE}) belongs to $C(0,1)\cap L^\infty(0,1)$ and satisfies $u(r)>0$ for all $r\in(0,1)$.
So condition (A0) is verified.

Now, by Theorem I of [\ref{IO}], we can obtain (i) and (ii). Finally, we prove (iii). Suppose now that (\ref{APE}) with
$\eta>\eta_1$ has a positive solution $v$, and let $u$ be a positive solution of (\ref{APE}) corresponding to $\eta_1(p)$. Lemma 3.1 and 3.2 imply that
every positive solution $w$ of $(\ref{APE})$ satisfies $w\in C^1[0,1]$ and $w'(1)<0$. By virtue of this fact and the fact that $tv$ is also
a solution of (\ref{APE}) for any real number $t$, we may assume without loss of generality that $u\leq v$. It is not difficult to verify
that $A$ and $B$ are monotone operators. The rest of proof is similar to that of [\ref{IO}, Theorem II]. \qed
\\ \\
\indent Let $\eta=\mu^{p-1}$, Lemma 3.3 shows the following result.\\
\\
\textbf{Lemma 3.4.} \emph{Put $\mu_1(p)=\left(\eta_1(p)\right)^{1/(p-1)}$.
Then we have that} \\

\emph{(i) (\ref{AP}) has no nontrivial solution for $\mu\in\left(0,\mu_1(p)\right)$;}

\emph{(ii) $\mu_1(p)$ is simple;}

\emph{(iii) (\ref{AP}) has a positive solution if and only if $\mu=\mu_1(p)$.}
\\ \\
\indent Moreover, we have the following result.
\\ \\
\textbf{Lemma 3.5.} \emph{If $(\mu,u)$ satisfies (\ref{AP}) with $\mu\neq \mu_1(p)$ and $u\not\equiv0$,
then $u$ must change sign.}\\
\\
\textbf{Proof.} Suppose that $u$ is not changing-sign. Without loss of
generality, we can assume that $u\geq 0$ in
$(0,1)$. Lemma 3.1 and 3.2 imply that $u>0$ in $(0,1)$. Lemma 3.4 implies $\mu=\mu_1(p)$ and $u=c v_1$
for some positive constant $c$, where $v_1$ is the positive eigenfunction
corresponding to $\mu_1(p)$ with $\left\Vert v_1\right\Vert=1$. This is a contradiction.\qed\\

\indent In addition, we also have that $\mu_1(p)$ is also isolated.
\\ \\
\textbf{Lemma 3.6.} \emph{$\mu_1(p)$ is the unique eigenvalue in $\left(0,\delta_p\right)$
for some $\delta_p>\mu_1(p)$.}
\\ \\
\textbf{Proof.} Lemma 3.4 has shown that $\mu_1(p)$ is left-isolated.
Assume by contradiction that there exists a sequence of eigenvalues
$\lambda_n\in\left(\mu_1(p), \delta_p\right)$ which converge to $\mu_1(p)$. Let $v_n$ be the
corresponding eigenfunctions.
Define
\begin{equation}
\psi_n:=\frac{v_n}{\left((p-1)\int_0^1 r^{p-2}a(r)\left\vert v_n\right\vert^p\,dr\right)^{1/p}}.\nonumber
\end{equation}
Clearly, $\psi_n$ are bounded in $W_c^{1,p}(0,1)$ so there exists a subsequence, denoted
again by $\psi_n$, and $\psi\in W_c^{1,p}(0,1)$ such that
$\psi_n\rightharpoonup \psi$ in $W_c^{1,p}(0,1)$ and $\psi_n\rightarrow \psi$ in $C^\alpha[0,1]$.
Since functional $f_1$ is sequentially weakly lower semi-continuous, we have that
\begin{equation}
\int_0^1 \left\vert \psi'\right\vert^p\,dr \leq\liminf_{n\rightarrow+\infty}\int_0^1
\left\vert \psi_n'\right\vert^p\,dr=\liminf_{n\rightarrow+\infty}\lambda_n^{p-1}=\mu_1^{p-1}(p).\nonumber
\end{equation}
On the other hand, $(p-1)\int_0^1 r^{p-2} a(r)\left\vert \psi_n\right\vert^p\,dr=1$ and
$\psi_n\rightarrow \psi$ in $C^\alpha[0,1]$ imply that
$(p-1)\int_0^1 r^{p-2}a(r)\vert \psi\vert^p\,dr=1$. Hence, $\int_0^1 \left\vert
\psi'\right\vert^p\,dr=\eta_1(p)$ via Lemma 3.3.
Then Lemma 3.1 and 3.3 show that $\psi>0$ in $(0,1)$. Thus $\psi_n\geq 0$
for $n$ large enough which contradicts the conclusion of Lemma 3.5.\qed
\\ \\
\indent Next, we show that the principle eigenvalue function $\mu_1:
[2,+\infty)\rightarrow \mathbb{R}$ is continuous.
\\ \\
\textbf{Lemma 3.7.} \emph{The eigenvalue function $\mu_1:[2,+\infty)\rightarrow \mathbb{R}$
is continuous.}
\\ \\
\textbf{Proof.} It is sufficient to show that $\eta_1(p):[2,+\infty)\rightarrow \mathbb{R}$
is continuous because of $\mu_1(p)=\left(\eta_1(p)\right)^{1/(p-1)}$.
From the variational characterization of $\eta_1(p)$ it follows that
\begin{equation}\label{cp1}
\eta_1(p)=\sup\left\{\lambda>0: \lambda(p-1)\int_0^1r^{p-2}a(r)\vert v\vert^p
\,dr\leq\int_0^1\left\vert v'\right\vert^p\,dr
\,\,\text{for all\,\,}v\in C_c^\infty[0,1]\right\},
\end{equation}
where $C_c^\infty[0,1]=\left\{v\in C^\infty[0,1]: v'(0)=v(1)=0\right\}$, as $C_c^\infty[0,1]$
is dense in $W_c^{1,p}(0,1)$ (see [\ref{A}]).

Let $\left\{p_j\right\}_{j=1}^\infty$ be a sequence in $[2, +\infty)$ which converge to $p\geq2$. We shall show
that
\begin{equation}\label{cp2}
\lim_{j\rightarrow+\infty}\eta_1\left(p_j\right)=\eta_1(p).
\end{equation}
To do this, let $v\in C_c^\infty[0,1]$. Then, due to (\ref{cp1}), we get that
\begin{equation}
\eta_1\left(p_j\right) \left(p_j-1\right)\int_0^1r^{p_j-2}a(r)\vert v\vert^{p_j}\,dr\leq\int_0^1\left\vert v'\right\vert^{p_j}\,dr.\nonumber
\end{equation}
On applying the Dominated Convergence Theorem we find that
\begin{equation}\label{cp3}
\limsup_{j\rightarrow+\infty}\eta_1\left(p_j\right) (p-1)\int_0^1r^{p-2}a(r)\vert
v\vert^{p}\,dr\leq\int_0^1\left\vert v'\right\vert^{p}\,dr.
\end{equation}
Relation (\ref{cp3}), the fact that $v$ is arbitrary and (\ref{cp1}) yield
\begin{equation}
\limsup_{j\rightarrow+\infty}\eta_1\left(p_j\right)\leq\eta_1(p).\nonumber
\end{equation}
\indent Thus, to prove (\ref{cp2}) it suffices to show that
\begin{equation}\label{cp4}
\liminf_{j\rightarrow+\infty}\eta_1\left(p_j\right)\geq\eta_1(p).
\end{equation}
Let $\left\{p_k\right\}_{k=1}^\infty$ be a subsequence of $\left\{p_j\right\}_{j=1}^\infty$ such that
$\underset{k\rightarrow+\infty}\lim\eta_1\left(p_k\right)=\underset{j\rightarrow+\infty}\liminf\eta_1\left(p_j\right)$.

Let us fix $\varepsilon_0>0$ so that $p-\varepsilon_0>1$ and for each
$0<\varepsilon<\varepsilon_0$ and $k\in \mathbb{N}$ large enough,
$p-\varepsilon<p_k<p+\varepsilon$. For $k\in \mathbb{N}$, let us choose $v_k\in W_c^{1,p_k}(0,1)$ such that $v_k>0$ in $(0,1)$,
\begin{equation}\label{cp5}
\int_0^1\left\vert v_k'\right\vert^{p_k}\,dr=1
\end{equation}
and
\begin{equation}\label{cp6}
\int_0^1\left\vert v_k'\right\vert^{p_k}\,dr=\eta_1\left(p_k\right)\left(p_k-1\right)\int_0^1 r^{p_k-2}a(r)\left\vert v_k\right\vert^{p_k}\,dr.
\end{equation}
(\ref{cp5}) shows that $\left\{v_k\right\}_{k=1}^\infty$ is a bounded sequence in $W_c^{1,p_k}(0,1)$, hence, in
$W_c^{1,p-\varepsilon}(0,1)$. Passing to a
subsequence if necessary, we can assume that $v_k \rightharpoonup v$ in
$W_c^{1,p-\varepsilon}(0,1)$ and hence
that $v_k \rightarrow v$ in $C^{\alpha}[0,1]$ with $\alpha=1-1/(p-\varepsilon)$
because the embedding of $W^{1,p-\varepsilon}(0,1)\hookrightarrow
C^{\alpha}[0,1]$ is compact. Thus,
\begin{equation}\label{cp8}
\left\vert v_k\right\vert^{p_k}\rightarrow\vert v\vert^{p}.
\end{equation}
\indent We note that (\ref{cp6}) implies that
\begin{equation}\label{cp9}
\eta_1\left(p_k\right)\left(p_k-1\right)\int_0^1 r^{p_k-2}a(r)\left\vert v_k\right\vert^{p_k}\,dr=1
\end{equation}
for all $k\in \mathbb{N}$. Thus letting $k\rightarrow+\infty$ in (\ref{cp9}) and using (\ref{cp8}), we find that
\begin{equation}\label{cp10}
\liminf_{j\rightarrow+\infty}\eta_1\left(p_j\right)(p-1)\int_0^1 r^{p-2}a(r)\vert v\vert^{p}\,dr=1.
\end{equation}
On the other hand, since $v_k\rightharpoonup v$ in $W_c^{1,p-\varepsilon}(0,1)$, from (\ref{cp5})
and the H\"{o}lder's inequality we obtain that
\begin{equation}
\left\Vert v'\right\Vert_{p-\varepsilon}^{p-\varepsilon}\leq\liminf_{k\rightarrow+\infty}\left\Vert v_k'
\right\Vert_{p-\varepsilon}^{p-\varepsilon}\leq1,\nonumber
\end{equation}
where $\Vert \cdot\Vert_p$ denotes the normal of $L^p(0,1)$.
Now, letting $\varepsilon\rightarrow 0^+$, we find
\begin{equation}\label{cp11}
\left\Vert v'\right\Vert_p\leq1.
\end{equation}
Clearly, (\ref{cp11}) and $v\in W_c^{1,p-\varepsilon}(0,1)$ follow
that $v\in W_c^{1,p}(0,1)$.

Consequently, combining (\ref{cp10}) with (\ref{cp11}) we obtain that
\begin{equation}
\liminf_{j\rightarrow+\infty}\eta_1\left(p_j\right)(p-1)\int_0^1 r^{p-2}a(r)\vert v\vert^{p}\,dr\geq
\int_0^1 \left\vert v'\right\vert^{p}\,dr.\nonumber
\end{equation}
This together with the variational characterization of $\eta_1(p)$ implies (\ref{cp4}) and hence
(\ref{cp2}). This concludes the proof of the lemma.\qed\\

We have shown that $\mu_1(p)$ is an isolated eigenvalue of (\ref{AP}). Moreover, we have the following result.
\\ \\
\textbf{Lemma 3.8.} For every interval $[a,  b]\subset [2,  +\infty)$ there is $\delta> 0$ such
that for all $p\in [a,  b]$ there is no eigenvalue of (\ref{AP}) in $\left(\mu_1(p),\mu_1(p)+\delta\right]$.
\\ \\
\textbf{Proof.} Suppose that the assertion of the proposition is not true. Then there are sequences
$\left\{ p_n\right\}_{n=1}^{+\infty}$ in $[2,+\infty)$, $\left\{\mu_n\right\}_{n=1}^{+\infty}$ in $\mathbb{R}^+$, and $\left\{v_n\right\}_{n=1}^{+\infty}$
in $X\setminus\{0\}$ such that $\lim_{n\rightarrow+\infty}p_n=p\in[2,+\infty)$, $\mu_n>\mu_1\left(p_n\right)$,
$\lim_{n\rightarrow+\infty}\left(\mu_n-\mu_1\left(p_n\right)\right)=0$, and
\begin{equation}
w_n=\int_1^r\varphi_{p_n'}\left(\int_s^0 \mu_n^{p_n-1}
\left(p_n-1\right)\tau^{p_n-2}a(\tau)\varphi_{p_n}\left(w_n\right)\,d\tau\right)\,ds,\,\, 0\leq r\leq1,\nonumber
\end{equation}
where $w_n=v_n/\left\Vert v_n\right\Vert$.
By Proposition 2.1 of [\ref{DPEM}], we can see that $\left\{w_n(r)\right\}$ is equicontinuous. By the Arzela-Ascoli theorem, we may assume that $w_n\rightarrow w$ in $X$. It follows that
\begin{equation}
w=\int_1^r\varphi_{p'}\left(\int_s^0 \mu_1(p)^{p-1}
(p-1)\tau^{p-2}a(\tau)\varphi_p(w)\,d\tau\right)\,ds,\,\, 0\leq r\leq1.\nonumber
\end{equation}
Then $\psi>0$ in $(0,1)$. Thus $v_n\geq 0$ for $n$ large enough which contradicts the assumption.\qed\\

We have known that $I-T^p_\mu$ is a completely continuous vector field in
$X$. Thus, the Leray-Schauder degree
$\deg\left(I-T^p_\mu, B_r(0),0\right)$ is well defined for
arbitrary $r$-ball $B_r(0)$ and $\mu\in\left(\mu_1(p),\mu_1(p)+\delta\right]$, where $\delta$ comes from Lemma 3.8.
\\ \\
\noindent\textbf{Theorem 3.1.} \emph{For fixed $p\geq 2$ and all $r>0$, we have that}
\begin{eqnarray}
\deg \left(I-T^p_\mu, B_r(0),0 \right)=\left\{
\begin{array}{lll}
1, \,\,\,\,\,\,&\text{if}\,\,\mu\in
\left(0,\mu_1(p)\right),\\
-1,\,\,&\text{if}\,\,\mu\in\left(\mu_1(p),\mu_1(p)+\delta\right).
\end{array}
\right.\nonumber
\end{eqnarray}
\noindent\textbf{Proof.} We shall only prove for the case $\mu>\mu_1(p)$
since the proofs of other cases are completely analogous. Assume that $\mu <\mu_1(p)+\delta$. Since the eigenvalue $\mu_1(p)$
depends continuously on $p$, there exist a continuous function
$\chi:[2,+\infty)\rightarrow\mathbb{R}$ and $q\geq2$ such that
$\mu_1(q) < \chi(q) < \mu_1(q)+ \delta$ and $\mu=\chi(p)$.
Set
\begin{equation}
d(q)=\deg\left(I-T_{\chi(q)}^q,B_r(0),0\right).\nonumber
\end{equation}
Since $T_\mu^2$ is compact and linear, by Theorem 8.10 of [\ref{De}], we have
\begin{equation}
d(2)=-1.\nonumber
\end{equation}
We shall show that $d(q)$ is locally constant in $[2,+\infty)$ and hence constant.
Define
\begin{equation}
G(q,v)=T_{\chi(q)}^qv.\nonumber
\end{equation}
Then by the Arzela-Ascoli theorem, we can show that $G:[2,+\infty)\times X\rightarrow X$ is completely continuous.
The invariance of the Leray-Shauder degree under a compact homotopy follows that $d(q)\equiv$\emph{constant} for $q\in\left[2,+\infty\right)$. So $\deg \left(I-T^p_\mu, B_r(0),0 \right)=d(p)=d(2)=-1$.\qed
\section{Global bifurcation result}

\quad\, With a simple transformation $v=-u$, problem (\ref{Mg}) can be written as
\begin{equation}\label{B1}
\left\{
\begin{array}{l}
\left(\left(-v'\right)^N\right)'=\lambda^NNr^{N-1}a(r)\left(v^N+g(v)\right),\,\, r\in(0,1),\\
v'(0)=v(1)=0.
\end{array}
\right.
\end{equation}
Let $X^+:=\{v\in X: v(r)\geq 0\}$ with the norm of $X$.
Define the map $T_g:X^+\rightarrow X^+$ by
\begin{equation}
T_gv(r)=\int_r^1\left(\int_0^sN\tau^{N-1}a(\tau)\left((v(\tau))^N+g(v(\tau))\right)
\,d\tau\right)^{1/N}\,ds,\,\, 0\leq r\leq1.\nonumber
\end{equation}
It is not difficult to verify that $T_g$ is continuous and compact.
Clearly, problem (\ref{B1}) can be equivalently written as
\begin{equation}
v=\lambda T_g v.\nonumber
\end{equation}
\indent Now, we show that the existence and uniqueness theorem is valid for problem (\ref{B1}).
\\ \\
\textbf{Lemma 4.1.} \emph{If $(\lambda, v)$ is a solution of (\ref{B1}) in $\mathbb{R}\times X^+$
and $v$ has a double zero, then $v \equiv 0$.}
\\ \\
\textbf{Proof.} Let $v$ be a solution of problem (\ref{B1}) and $r_*\in[0, 1]$ be a double zero.
We note that
\begin{equation}
v(r)=\lambda\int_r^{r_*}\left(\int_{r_*}^sN\tau^{N-1}a(\tau)\left((v(\tau))^N+g(v(\tau))\right)
\,d\tau\right)^{1/N}\,ds.\nonumber
\end{equation}
Firstly, we consider $r\in[0, r_*]$. Then we have that
\begin{eqnarray}
\vert v(r)\vert&\leq&\lambda\left(\int_{r}^{r_*}N\tau^{N-1}a(\tau)\left\vert\left((v(\tau))^N
+g(v(\tau))\right)\right\vert\,d\tau\right)^{1/N},\nonumber
\end{eqnarray}
furthermore,
\begin{eqnarray}
\vert v(r)\vert^N&\leq&\lambda^N\int_{r}^{r_*}N\tau^{N-1}a(\tau)\left\vert\left((v(\tau))^N
+g(v(\tau))\right)\right\vert\,d\tau \nonumber\\
&\leq&\lambda^N\int_r^{r^*}
 N \tau^{N-1}a(\tau)\left\vert 1+\frac{g(v(\tau))}{(v(\tau))^N}\right\vert \vert v(\tau)\vert^N\,d\tau.\nonumber
\end{eqnarray}
According to the assumptions on $g$, for any $\varepsilon>0$, there exists a constant $\delta>0$ such that
\begin{equation}
\vert g(s)\vert\leq \varepsilon s^N\nonumber
\end{equation}
for any $s\in[0,\delta]$.
Hence, we have that
\begin{equation}
\vert v(r)\vert^N\leq \lambda^N\int_r^{r^*}
N a^0\left(1+\varepsilon+\max_{s\in\left[\delta,\Vert v\Vert\right]}
\left\vert\frac{g(s)}{s^N}\right\vert\right) \vert v(\tau)\vert^N\,d\tau,\nonumber
\end{equation}
where $a^0=\max_{r\in[0,1]}a(r)$.
By the modification of the Gronwall-Bellman inequality [\ref{ILL}, Lemma 3.2], we get $v \equiv 0$ on $\left[0, r^*\right]$.
Similarly, using the Gronwall-Bellman inequality [\ref{Bre}, \ref{E}], we can get $v \equiv 0$ on $\left[r^*, 1\right]$
and the proof is completed.\qed\\

Now, we consider the following eigenvalue problem
\begin{equation}\label{B2}
\left\{
\begin{array}{l}
\left(\left(-v'\right)^N\right)'=\lambda^NNr^{N-1}a(r)v^N,\,\, r\in(0,1),\\
v'(0)=v(1)=0.
\end{array}
\right.
\end{equation}
The same proof as in Theorem 1.1 of [\ref{L1}], we can show that problem (\ref{B2}) possesses the first eigenvalue $\lambda_1$ which is positive,
simple and the corresponding eigenfunctions are positive in (0,1) and concave on $[0,1]$. Moreover, we also have the following result.\\
\\
\textbf{Lemma 4.2.} \emph{If $(\mu,\varphi)\in(0,+\infty)\times \left(C^2[0,1]
\setminus\{0\}\right)$ satisfies (\ref{B2}) and $\mu\neq\lambda_1$,
then $\varphi$ must change sign.}\\
\\
\textbf{Proof.} By way of contradiction, we may suppose that $\varphi$ is not changing-sign. Without loss of
generality, we can assume that $\varphi\geq 0$ in
$(0,1)$. Lemma 4.1 follows that $\varphi>0$ in $(0,1)$. It implies $\mu=\lambda_1$ and $\varphi=\theta \psi_1$
for some positive constant $\theta$, where $\psi_1$ is the positive eigenfunction
corresponding to $\lambda_1$ with $\Vert \psi_1\Vert=1$. We have a contradiction.\qed\\
\\
\indent Next, we show that $\lambda_1$ is also isolated.
\\ \\
\textbf{Lemma 4.3.} \emph{$\lambda_1$ is isolated; that is to say, $\lambda_1$ is
the unique eigenvalue in $(0,\delta)$ for some $\delta>\lambda_1$.}\\ \\
\textbf{Proof.} We have known that $\lambda_1$ is left-isolated.
Assume by contradiction that there exists a sequence of eigenvalues
$\lambda_n\in\left(\lambda_1, \delta\right)$ which converge to $\lambda_1$. Let $v_n$ be the
corresponding eigenfunctions. Let $w_n:=v_n/\left\Vert v_n\right\Vert_{C^1[0,1]}$,
then $w_n$ should be the solutions of the following problem
\begin{equation}
w=\lambda_n\int_r^1\left(\int_0^sN\tau^{N-1}a(\tau)w^N\,d\tau\right)^{1/N}\,ds.\nonumber
\end{equation}
Clearly, $w_n$ are bounded in $C^1[0,1]$ so there exists a subsequence, denoted again by $w_n$,
and $\psi\in X$ such that
$w_n\rightarrow\psi$ in $X$.
It follows that
\begin{equation}
\psi=\lambda_1\int_r^1\left(\int_0^sN\tau^{N-1}a(\tau)\psi^N\,d\tau\right)^{1/N}\,ds.\nonumber
\end{equation}
It follows that $\psi=\theta\psi_1$ for some positive constant $\theta$
in $(0,1)$. Thus $w_n\geq 0$ for $n$ large enough contradicts $v_n$
changing-sign in $(0,1)$ which is implied by Lemma 4.2.\qed\\

Define $T_N:X^+\rightarrow X^+$ by
\begin{equation}
T_Nv:=\int_r^1\left(\int_0^sN\tau^{N-1}a(\tau)v^N\,d\tau\right)^{1/N}\,ds,\,\, 0\leq r\leq1.\nonumber
\end{equation}
Clearly, $I-T_N$ is a completely continuous vector field in
$X^+$. Thus, the Leray-Schauder degree
$\deg\left(I-T_N, B_r(0),0\right)$ is well defined for
arbitrary $r$-ball $B_r(0)$ of $X^+$ and $\mu\in (0,\delta)\setminus\{\lambda_1\}$, where $\delta$ comes from Lemma 4.3.
\\ \\
\noindent\textbf{Lemma 4.4.} \emph{Let $\lambda$ be a constant with
$\lambda\in(0,\delta)$. Then for
arbitrary $r>0$,}
\begin{eqnarray}
\deg \left(I-\lambda T_N, B_r(0),0 \right)=\left\{
\begin{array}{lll}
1, \,\,\,\,\,\,&\text{if}\,\,\lambda\in
\left(0,\lambda_1\right),\\
-1,\,\,&\text{if}\,\,\lambda\in\left(\lambda_1,\delta\right).
\end{array}
\right.\nonumber
\end{eqnarray}
\noindent\textbf{Proof.} Taking $p=N+1$ and $\mu=\lambda$ in $T_\mu^p$, we can see that
$\lambda_1=\mu_1(p)$. Furthermore, it is
not difficult to verify that $\lambda T_N(v)=T_\mu^p(v)$ for any $v\in X^+$.
By Theorem 3.1, we can deduce this lemma.\qed
\\ \\
\textbf{Theorem 4.1.} \emph{$\left(\lambda_1,0\right)$ is a bifurcation
point of (\ref{B1}) and the associated bifurcation branch $\mathcal{C}$ in $\mathbb{R}\times X^+$
whose closure contains $\left(\lambda_1, 0\right)$ is either unbounded or contains a pair $\left(\overline{\lambda}, 0\right)$
where $\overline{\lambda}$ is an eigenvalue of (\ref{B2}) and $\overline{\lambda}\neq \lambda_1$.}
\\ \\
\textbf{Proof.} Suppose that $\left(\lambda_1, 0\right)$
is not a bifurcation point of problem (\ref{B1}). Then
there exist $\varepsilon > 0$, $\rho_0> 0$ such that for $\left\vert \lambda-\lambda_1\right\vert\leq\varepsilon$
and $0<\rho < \rho_0$ there is no
nontrivial solution of the equation
\begin{equation}
v-\lambda T_gv=0\nonumber
\end{equation}
with $\Vert v\Vert=\rho$. From the invariance of the degree under a compact
homotopy we obtain that
\begin{equation}\label{edc}
\text{deg}\left(I-\lambda T_g,B_\rho(0),0\right)\equiv constant
\end{equation}
for $\lambda\in\left[\lambda_1-\varepsilon,\lambda_1+\varepsilon\right]$.

By taking $\varepsilon$ smaller if necessary, in view of Lemma 4.3, we can assume
that there is no eigenvalue
of (\ref{B2}) in $\left(\lambda_1,\lambda_1+\varepsilon\right]$. Fix $\lambda\in\left(\lambda_1,\lambda_1+\varepsilon\right]$.
We claim that the equation
\begin{equation}\label{es}
v-\lambda\int_r^1\left(\int_0^sN\tau^{N-1}a(\tau)\left(v^N+tg(v)\right)\,d\tau\right)^{1/N}\,ds=0
\end{equation}
has no solution $v$ with $\Vert v\Vert=\rho$ for every $t\in[0, 1]$ and $\rho$
sufficiently small.
Suppose on the contrary, let $\left\{v_n\right\}$ be the nontrivial solutions of (\ref{es})
with $\left\Vert v_n\right\Vert\rightarrow 0$
as $n\rightarrow+\infty$.

Let $w_n:=v_n/\left\Vert v_n\right\Vert$, then $w_n$ should be the solutions of the following problem
\begin{equation}\label{B3}
w(t)=\lambda\int_r^1\left(\int_0^sN\tau^{N-1}a(\tau)\left(w^N+t\frac{g(v)}
{\left\Vert v_n\right\Vert^N}\right)\,d\tau\right)^{1/N}\,ds.
\end{equation}
Let
\begin{equation}
\widetilde{g}(v)=\max_{0\leq s\leq v}\vert g(s)\vert,\nonumber
\end{equation}
then $\widetilde{g}$ is nondecreasing with respect to $v$ and
\begin{equation}\label{eg0+}
\lim_{ v\rightarrow 0^+}\frac{\widetilde{g}(v)}{
v^{N}}=0.
\end{equation}
Further it follows from (\ref{eg0+}) that
\begin{equation}\label{egn0}
\frac{\vert g(v)\vert}{\Vert v\Vert^{N}} \leq \frac{\widetilde{ g}(v)}
{\Vert v\Vert^{N}}\leq \frac{
\widetilde{g}(\Vert v\Vert)}{\Vert v\Vert^{N}}\rightarrow0\,\,  \text{as}\,\, \Vert
v\Vert\rightarrow 0.
\end{equation}
By (\ref{B3}), (\ref{egn0}) and the compactness of $T_g$, we obtain that for
some convenient subsequence
$w_n\rightarrow w_0$ as $n\rightarrow+\infty$. Now $\left(\lambda,w_0\right)$ verifies
problem (\ref{B2}) and $\left\Vert w_0\right\Vert = 1$. This implies that $\lambda$ is an eigenvalue of (\ref{B2}).
This is a contradiction.

From the invariance of the degree under
homotopies and Lemma 4.4 we then obtain that
\begin{equation}\label{edFk}
\deg\left(I-\lambda T_g(\cdot), B_r(0),0\right)=
\deg\left(I-\lambda T_N(\cdot), B_r(0),0\right)=-1.
\end{equation}
Similarly, for $\lambda\in \left[\lambda_1 - \varepsilon, \lambda_1\right)$ we find that
\begin{equation}\label{edFk1}
\deg\left(I-\lambda T_g(\cdot), B_r(0),0\right)=1.
\end{equation}
Relations (\ref{edFk}) and (\ref{edFk1}) contradict (\ref{edc}) and hence $\left(\lambda_1, 0\right)$ is a
bifurcation point of problem (\ref{B1}).

By standard arguments in global bifurcation theory (see [\ref{R}]), we can
show the existence of
a global branch of solutions of problem (\ref{B1}) emanating from
$\left(\lambda_1, 0\right)$. Our conclusion is proved.\qed\\

Next, we shall prove that the first choice of the alternative of Theorem 4.1 is
the only possibility. Let $P^+$ denote the set of functions
in $X^+$ which are positive in (0,1). Set $K^{+}=\mathbb{R}\times P^{+}$ under the product topology.
\\ \\
\textbf{Theorem 4.2.} \emph{There exists an unbounded
continuum $\mathcal{C}\subseteq \left(K^+\cup\{\left(\lambda_1,0\right)\}\right)$ of solutions to problem (\ref{B1})
emanating from $\left(\lambda_1,0\right)$. Moreover, $\left(\lambda_1,0\right)$ is the unique bifurcation point
from $(\lambda,0)$ of the positive solutions of problem (\ref{B1}).}\\ \\
\textbf{Proof.} For any $(\lambda,v)\in \mathcal{C}$, Lemma 4.1 implies that either $v\equiv 0$ or $v>0$ in (0,1).
Thus, we have $\mathcal{C}\subseteq \left(K^+\cup\{\left(\lambda_1,0\right)\}\right)$.

Now, we prove that the first choice of the alternative of Theorem 4.1 is
the only possibility. Suppose on the contrary, if there exists $\left(\lambda_n,v_n\right)
\rightarrow\left(\overline{\lambda},0\right)$
when $n\rightarrow+\infty$ with $\left(\lambda_n,v_n\right)\in \mathcal{C}$,
$v_n \not\equiv 0$ and $\overline{\lambda}$
is another eigenvalue of (\ref{B2}).
Let $w_n :=v_n/\left\Vert v_n\right\Vert$, then $w_n$ should be the solutions of the following problem
\begin{equation}\label{evs}
w=\lambda_n\int_r^1\left(\int_0^sN\tau^{N-1}a(\tau)\left(w^N+\frac{g(v)}
{\left\Vert v_n\right\Vert^N}\right)\,d\tau\right)^{1/N}\,ds.
\end{equation}
By an argument similar to that of Theorem 4.1, we obtain that for some convenient
subsequence $w_n\rightarrow w_0$ as $n\rightarrow+\infty$. It is easy to
see that $\left(\overline{\lambda}, w_0\right)$ verifies problem (\ref{B2})
and $\left\Vert w_0\right\Vert = 1$. Lemma 4.2 follows $w_0$ must change sign, and as a
consequence for some $n$ large enough, $w_n$ must change sign, and
this is a contradiction.\qed
\\

Finally, we give a key lemma that will be used later.\\
\\
\textbf{Lemma 4.5.} \emph{Let $b_2(r)\geq b_1(r)>0$ for $r\in(0,1)$ and $b_i(r)\in C([0,1])$, $i=1,2$. Also let $u_1$, $u_2$
be solutions of the following differential problems
\begin{equation}\label{411}
\left\{
\begin{array}{l}
\left(\left(-u'\right)^N\right)'=b_i(r)u^N,\,\, i=1,2,\\
u'(0)=u(1)=0,
\end{array}
\right.
\end{equation}
respectively. If $u_1(r)\neq 0$ in $(0,1)$, then either there exists $\tau\in (0,1)$
such that $u_2(\tau)=0$ or $b_2=b_1$ and $u_2(r)=\mu u_1(r)$ for some constant $\mu\neq0$ and a.e. $r\in(0,1)$.}
\\ \\
\textbf{Proof.} If $u_2(r)\neq0$ in $(0,1)$, then we can assume without loss of generality that
$u_1(r)>0$, $u_2(r)>0$ in $(0,1)$. Then from problem (3.11), we can easily show that $u_1$ and $u_2$ are strictly decreasing concave functions in $(0,1)$.
Moreover, it is easy to check that the conclusion of Lemma 4.1 is also valid for problem (3.11). So we have $u_1'(r)<0$ and $u_2'(r)<0$ for $r\in(0,1]$.

By some simple calculations, we have that
\begin{eqnarray}
&&\int_0^1\left(\frac{u_1^{N+1}\left(-u_2'\right)^N}{u_2^N}-u_1\left(-u_1'\right)^N\right)'dr=\nonumber\\
&&\int_0^1\left(wu_1^{N+1}+\left(\left(-u_1'\right)^{N+1}+N\left(\frac{-u_1u_2'}{u_2}\right)^{N+1}-(N+1)u_1^Nu_1'\left(\frac{-u_2'}{u_2}\right)^N\right)\right)dr,
\end{eqnarray}
where $w=b_2-b_1$.
The left-hand side of (3.12) equals
\begin{equation}
\lim_{r\rightarrow 1^-}\frac{u_1^{N+1}\left(-u_2'\right)^N}{u_2^N}:=H.\nonumber
\end{equation}
We prove that $H=0$. By the L'Hospital rule, we have that
\begin{eqnarray}
H&=&\lim_{r\rightarrow 1^-}\frac{u_1^{N+1}\left(-u_2'\right)^N}{u_2^N}
=\lim_{r\rightarrow 1^-}\frac{(N+1)u_1^Nu_1'\left(-u_2'\right)^N+u_1^{N+1}\left(\left(-u_2'\right)^N\right)'}
{Nu_2^{N-1}u_2'}\nonumber\\
&=&\lim_{r\rightarrow 1^-}\frac{(N+1)u_1^Nu_1'\left(-u_2'\right)^N+u_1^{N+1}b_2u_2^N}
{Nu_2^{N-1}u_2'}\nonumber\\
&=&\lim_{r\rightarrow 1^-}\frac{(N+1)u_1^Nu_1'\left(-u_2'\right)^N}
{Nu_2^{N-1}u_2'}+\lim_{r\rightarrow 1^-}\frac{u_1^{N+1}b_2u_2^N}
{Nu_2^{N-1}u_2'}\nonumber\\
&=&\lim_{r\rightarrow 1^-}\frac{(N+1)u_1'\left(-u_2'\right)^N}
{Nu_2'}\lim_{r\rightarrow 1^-}\frac{u_1^N}
{u_2^{N-1}}.\nonumber
\end{eqnarray}
If $N=1$, then $H=0$. If $1<N\leq 2$, applying the L'Hospital rule again, we obtain that
\begin{equation}
\lim_{r\rightarrow 1^-}\frac{u_1^{N}}{u_2^{N-1}}=\lim_{r\rightarrow 1^-}\frac{Nu_1'}{(N-1)u_2'}\lim_{r\rightarrow 1^-}\frac{u_1^{N-1}}{u_2^{N-2}}.\nonumber
\end{equation}
This implies that $H=0$. If $k-1<N\leq k$, then we continue this process $k$ times to obtain $H=0$.

Therefore, the left-hand side of (3.12) equals zero. Hence
the right-hand side of (3.12) also equals zero.
The Young's inequality implies that
\begin{equation}
\left(-u_1'\right)^{N+1}+N\left(\frac{-u_1u_2'}{u_2}\right)^{N+1}-(N+1)u_1^Nu_1'\left(\frac{-u_2'}{u_2}\right)^N\geq 0,\nonumber
\end{equation}
and the equality holds if and only if
\begin{equation}
\left(\frac{-u_1'}{u_1}\right)^{N+1}=\left(\frac{-u_2'}{u_2}\right)^{N+1}.\nonumber
\end{equation}
It follows that there exists a constant $\mu\neq0$ such that $u_2=\mu u_1$ and $b_2=b_1$.\qed\\

\indent As an immediate consequence, we obtain the following Sturm type comparison lemma.
\\ \\
\textbf{Lemma 4.6.} \emph{Let $b_i(r)\in C([0,1])$, $i=1,2$ such that $b_2(r)\geq b_1(r)>0$ for $r\in(0,1)$ and the inequality is strict
on some subset of positive measure in $(0,1)$. Also let $u_1$, $u_2$
are solutions of (\ref{411}) with $i=1,2$,
respectively. If $u_1\neq 0$ in $(0,1)$, then $u_2$ has at least one zero in $(0,1)$.}

\section{Convex solutions}

\quad\, In this section, we shall investigate the existence and multiplicity of convex
solutions of problem (\ref{MO}). With a simple transformation $v=-u$, problem (\ref{MO})
can be written as
\begin{equation}\label{C1}
\left\{
\begin{array}{l}
\left(\left(-v'(r)\right)^N\right)'=\lambda^NNr^{N-1}a(r)f(v(r)),\,\, r\in(0,1),\\
v'(0)=v(1)=0.
\end{array}
\right.
\end{equation}
Define the map $T_f:X^+\rightarrow X^+$ by
\begin{equation}
T_fv(r)=\int_r^1\left(\int_0^sN\tau^{N-1}a(\tau)f(v(r))\,d\tau\right)^{1/N}
\,ds,\,\, 0\leq r\leq1.\nonumber
\end{equation}
Similar to $T_g$, $T_f$ is continuous and compact. Clearly, problem (\ref{C1})
can be equivalently written as
\begin{equation}
v=\lambda T_f v.\nonumber
\end{equation}
\indent Let $f_0, f_\infty\in \mathbb{R}\setminus \mathbb{R}^-$ be such that
\begin{equation}
f_0=\lim_{ s\rightarrow0^+}\frac{f(s)}{s^N}\,\,\text{and}\,\,f_\infty
=\lim_{s\rightarrow+\infty}\frac{f(s)}{s^N}.\nonumber
\end{equation}
Throughout this section, we always suppose that $f$ satisfies the following signum condition\\

(f1) \emph{$f\in C(\mathbb{R}\setminus \mathbb{R}^-,\mathbb{R}\setminus \mathbb{R}^-)$ with
$f(s)s^N>0$ for $s>0$.}\\

Applying Theorem 4.2, we shall establish the existence of convex solutions of (\ref{MO}) as follows.\\ \\
\textbf{Theorem 5.1.} \emph{If $f_0\in(0,+\infty)$ and $f_\infty\in(0,+\infty)$,
then for any $\lambda\in\left(\lambda_1/f_\infty,\lambda_1/f_0\right)$
or $\lambda\in\left(\lambda_1/f_0,\lambda_1/f_\infty\right)$, (\ref{MO})
has at least one solution $u$ such that it is negative,
strictly convex in $(0,1)$.}
\\ \\
\textbf{Proof.} It suffices to prove that (\ref{C1}) has at least one solution $v$ such
that it is positive, strictly concave in $(0,1)$.

Clearly, $f_0\in(0,+\infty)$ implies $f(0)=0$. Hence, $v=0$ is always the solution of problem (\ref{C1}). Let $\zeta\in C\left(\mathbb{R}\setminus\mathbb{R}^-,\mathbb{R}\setminus\mathbb{R}^-\right)$ be such that $f(s)=f_0 s^N+\zeta(s)$
with $\lim_{s\rightarrow0^+}\zeta(s)/s^N=0.$
Applying Theorem 4.2 to (\ref{C1}), we have that
there exists an unbounded continuum $\mathcal{C}$ emanating
from $\left(\lambda_1/f_0, 0\right)$, such that
\begin{equation}
\mathcal{C}\subseteq \left(\left\{\left(\lambda_1,0\right)\right\}
\cup\left(\mathbb{R}\times P^+\right)\right).\nonumber
\end{equation}
\indent To complete this theorem, it will be enough to show that $\mathcal{C}$ joins
$\left(\lambda_1/f_0, 0\right)$ to
$\left(\lambda_1/f_\infty, +\infty\right)$. Let
$\left(\mu_n, v_n\right) \in \mathcal{C}$ satisfy
$\mu_n+\left\Vert v_n\right\Vert\rightarrow+\infty.$
We note that $\mu_n >0$ for all $n \in \mathbb{N}$ since (0,0) is
the only solution of (\ref{C1}) for $\lambda= 0$ and
$\mathcal{C}\cap\left(\{0\}\times X^+\right)=\emptyset$.

We divide the rest of proofs into two steps.

\emph{Step 1.} We show that there exists a constant $M$ such that $\mu_n
\in(0,M]$ for $n\in \mathbb{N}$ large enough.

On the contrary, we suppose that $\lim_{n\rightarrow +\infty}\mu_n=+\infty.$
On the other hand, we note that
\begin{equation}
\left(\left(-v_n'(r)\right)^N\right)'=\mu_n^N N r^{N-1}
\frac{f(v_n)}{v_n^N}v_n^N.\nonumber
\end{equation}
The signum condition (f1) combining $f_0\in(0,+\infty)$ and $f_\infty\in(0,+\infty)$ implies that there exists a positive
constant $\varrho$ such that $f(v_n)/v_n^N\geq\varrho$ for
any $r\in(0,1)$.
By Lemma 4.6, we get $v_n$ must change sign in $(0,1)$ for $n$
large enough, and this contradicts the fact that $v_n\in \mathcal{C}$.

\emph{Step 2.} We show that $\mathcal{C}$ joins
$\left(\lambda_1/f_0, 0\right)$ to
$\left(\lambda_1/ f_\infty, +\infty\right)$.

It follows from \emph{Step 1} that $\left\Vert v_n\right\Vert\rightarrow+\infty.$
Let $\xi\in C\left(\mathbb{R}\setminus\mathbb{R}^-,\mathbb{R}\setminus\mathbb{R}^-\right)$ be such that $f(s)=f_\infty s^N+\xi(s).$
Then $\lim_{s\rightarrow+\infty}\xi(s)/s^N=0.$
Let $\widetilde{\xi}(v)=\max_{0\leq s\leq v}\vert \xi(s)\vert.$
Then $\widetilde{\xi}$ is nondecreasing.
Set $\overline{\xi}(v)=\max_{v/2\leq s\leq v}\vert \xi(s)\vert.$ Then we have
\begin{equation}
\lim_{v\rightarrow +\infty}\frac{\overline{\xi}(v)}{v^{N}}=0\,\,\text{and}\,\,\widetilde{\xi}(v)\leq\widetilde{\xi}\left(\frac{v}{2}\right)+\overline{\xi}(v).\nonumber
\end{equation}
It follows that
\begin{equation}\label{eu0}
\lim_{v\rightarrow +\infty}\frac{\widetilde{\xi}(v)}{v^{N}}=0.
\end{equation}
\indent We divide the equation
\begin{equation}
\left(\left(-v_n'\right)^N\right)'-\mu_n^Nf_\infty r^{N-1}a(r)v_n^N
=\mu_n^N r^{N-1}a(r)\xi(v_n)\nonumber
\end{equation}
by $\left\Vert v_n\right\Vert^N$ and set $\overline{v}_n = v_n/\left\Vert v_n\right\Vert$.
Since $\overline{v}_n$ are bounded in $X^+$,
after taking a subsequence if
necessary, we have that $\overline{v}_n \rightharpoonup \overline{v}$
for some $\overline{v} \in X^+$. Moreover, from
(\ref{eu0}) and the fact that $\widetilde{\xi}$ is nondecreasing,
we have that
\begin{equation}\label{C2}
\lim_{n\rightarrow+\infty}\frac{ \xi\left(v_n(r)\right)}{\left\Vert v_n\right\Vert^{N}}=0
\end{equation}
since
\begin{equation}
\frac{ \vert\xi\left(v_n(r)\right)\vert}{\left\Vert v_n\right\Vert^{N}}\leq\frac{ \widetilde{\xi}
(\left\vert v_n(r)\right\vert)}{\left\Vert v_n\right\Vert^{N}}
\leq\frac{ \widetilde{\xi}(\left\Vert v_n(r)\right\Vert)}{\left\Vert v_n\right\Vert^{N}}.\nonumber
\end{equation}
\indent By the continuity and compactness of $T_f$, it follows that
\begin{equation}
\left(\left(-\overline{v}'\right)^N\right)'-\overline{\lambda}^Nf_\infty r^{N-1}a(r)
\overline{v}^N=0,\nonumber
\end{equation}
where
$\overline{\lambda}=\underset{n\rightarrow+\infty}\lim\lambda_n$, again
choosing a subsequence and relabeling it if necessary.

It is clear that $\left\Vert \overline{v}\right\Vert=1$ and $\overline{v}\in
\overline{\mathcal{C}}\subseteq
\mathcal{C}$ since $\mathcal{C}$ is closed in $\mathbb{R}\times X^+$.
Therefore, $\overline{\lambda}f_\infty=\lambda_1$, so $\overline{\lambda}=\lambda_1/ f_\infty.$
Therefore, $\mathcal{C}$ joins $\left(\lambda_1/
f_0, 0\right)$ to $\left(\lambda_1/f_\infty,
+\infty\right)$.\qed\\
\\
\textbf{Remark 5.1.} From the proof of Theorem 5.1, we can see
that if $f_0, f_\infty\in(0,+\infty)$ then there exist $\lambda_2>0$ and $\lambda_3>0$ such that
(\ref{MO}) has at least one strictly convex solution for all $\lambda\in\left(\lambda_2,\lambda_3\right)$
and has no convex solution for all
$\lambda\in\left(0,\lambda_2\right)\cup\left(\lambda_3,+\infty\right)$.
\\ \\
\textbf{Proof.} Clearly, $f_0, f_\infty\in(0,+\infty)$ implies that there exists a positive constant $M$ such that
\begin{equation}
\left\vert\frac{f(s)}{s^N}\right\vert\leq M \,\,\text{for any}\,\, s>0.\nonumber
\end{equation}
It is sufficient to show that there exists $\lambda_2>0$ such that
(\ref{MO}) has no convex solution for all
$\lambda\in\left(0,\lambda_2\right)$. Suppose on the contrary that there exists one pair $\left(\mu, v\right)\in\mathcal{C}$ such that
$\mu\in\left(0,1/\left(M^{1/N }a^0\right)\right)$.
Let $w =v/\left\Vert v\right\Vert$. Obviously, one has that
\begin{eqnarray}
1=\left\Vert w\right\Vert=\left\Vert\mu\int_r^1\left(\int_0^sN\tau^{N-1}a(\tau)\left(\frac{f\left(v\right)}
{\left\Vert v\right\Vert^N}\right)\,d\tau\right)^{1/N}\,ds\right\Vert
\leq M^{1/N}a^0\mu<1. \nonumber
\end{eqnarray}
This is a contradiction.\qed\\

From the proof of Theorem 5.1 and  Remark 5.1, we can deduce the following two corollaries.
\\
\\
\textbf{Corollary 5.1.} \emph{Assume that there exists a positive constant $\rho>0$ such that
\begin{equation}
\frac{f(s)}{s^N}\geq \rho\nonumber
\end{equation}
for any $s>0$.
Then there exists $\zeta_*>0$ such that problem (\ref{MO}) has no convex solution for
any $\lambda\in\left(\zeta_*,+\infty\right)$.}
\\ \\
\textbf{Corollary 5.2.} \emph{Assume that there exists a positive constant $\varrho>0$ such that
\begin{equation}
\left\vert\frac{f(s)}{s^N}\right\vert\leq\varrho\nonumber
\end{equation}
for any $s>0$.
Then there exists $\eta_*>0$ such that problem (\ref{MO}) has no convex solution for
any $\lambda\in\left(0,\eta_*\right)$.}
\\ \\
\textbf{Theorem 5.2.} \emph{If $f_0\in(0,+\infty)$ and $f_\infty=0$,
then for any $\lambda\in\left(\lambda_1/f_0,+\infty\right)$, (\ref{MO})
has at least one solution $u$ such that it is negative,
strictly convex in $(0,1)$.}
\\ \\
\textbf{Proof.} In view of Theorem 5.1, we only need to show that
$\mathcal{C}$ joins
$\left(\lambda_1/f_0, 0\right)$ to
$\left(+\infty, +\infty\right)$.

We first show that $\mathcal{C}$ is unbounded in the direction of $X^+$.
Suppose on the contrary that there exists $M>0$ such that $\Vert v\Vert\leq M$ for any $(\lambda,v)\in \mathcal{C}$.
Let $w =v/\left\Vert v\right\Vert$. Obviously, one has that
\begin{eqnarray}
1\geq w&=&\lambda\int_r^1\left(\int_0^sN\tau^{N-1}a(\tau)\left(\frac{f\left(v\right)}
{\left\Vert v\right\Vert^N}\right)\,d\tau\right)^{1/N}\,ds\nonumber \\
&=&\lambda\int_r^1\left(\int_0^sN\tau^{N-1}a(\tau)\left(\frac{f\left(v\right)}
{v^N}\right)\left(\frac{v^N}
{\left\Vert v\right\Vert^N}\right)\,d\tau\right)^{1/N}\,ds. \nonumber
\end{eqnarray}
The signum condition (f1) combining $f_0\in(0,+\infty)$ and $\Vert v\Vert\leq M$ implies that there exists a positive
constant $\rho$ such that $f(v)/v^N\geq\rho$ for
any $r\in(0,1)$. For any $\alpha\in(0,1)$, Lemma 2.2 of [\ref{HW}] shows that
\begin{eqnarray}
\min_{r\in[\alpha,1-\alpha]} v(r)\geq \alpha \Vert v\Vert. \nonumber
\end{eqnarray}
So for any $r\in[\alpha,1-\alpha]$, we have
\begin{eqnarray}
1\geq w(r)=\lambda \alpha \rho^{1/N}N^{1/N}\xi,\nonumber
\end{eqnarray}
where $\xi=\min_{r\in[\alpha,1-\alpha]}\int_r^1\left(\int_0^s\tau^{N-1}a(\tau)\,d\tau\right)^{\frac{1}{N}}\,ds>0$.
Hence, we have
\begin{eqnarray}
\lambda\leq \frac{1}{\alpha \rho^{1/N}N^{1/N}\xi}.\nonumber
\end{eqnarray}
This is a contradiction.

Next, we show that the unique  blow up point (see Definition 1.1 of [\ref{S}]) is $\left(+\infty,0\right)$. Suppose on the contrary that there exists $\mu_M$ such that $\left(\mu_M,0\right)$ is a blow up point and $\mu_M<+\infty$.
Then there exists a sequence $\left\{\mu_n, v_n\right\}$ such that
$\underset{n\rightarrow +\infty}{\lim}\mu_n=\mu_{M}$ and
$\underset{n\rightarrow +\infty}{\lim}\left\Vert v_n\right\Vert=+\infty$ as
$n\rightarrow+\infty$. Let $w_n =v_n/\left\Vert v_n\right\Vert$ and $w_n$
should be the solutions of the following problem
\begin{equation}
w=\mu_n\int_r^1\left(\int_0^sN\tau^{N-1}a(\tau)\left(\frac{f\left(v_n\right)}
{\left\Vert v_n\right\Vert^N}\right)\,d\tau\right)^{1/N}\,ds.\nonumber
\end{equation}
Similar to (\ref{C2}), we can show that
\begin{equation}
\lim_{n\rightarrow+\infty}\frac{ f\left(v_n(r)\right)}{\left\Vert v_n\right\Vert^{N}}=0.\nonumber
\end{equation}
By the compactness of $T_f$, we obtain that for some convenient subsequence
$w_n\rightarrow w_0$ as $n\rightarrow+\infty$. Letting $n\rightarrow+\infty$,
we obtain that
$w_0\equiv0$. This contradicts $\left\Vert w_0\right\Vert=1$. \qed\\
\\
\textbf{Remark 5.2.} Under the assumptions of Theorem 5.2, in view of Corollary 5.2, we can see that there exists $\lambda_4>0$ such that
problem (\ref{MO}) has at least one strictly convex solution for all $\lambda\in\left(\lambda_4,+\infty\right)$
and has no convex solution for all $\lambda\in\left(0,\lambda_4\right)$.
\\ \\
\textbf{Theorem 5.3.} \emph{If $f_0\in(0,+\infty)$ and $f_\infty=\infty$,
then for any $\lambda\in\left(0,\lambda_1/f_0\right)$, (\ref{MO}) has at least one
solution $u$ such that it is negative, strictly convex in $(0,1)$.}
\\ \\
\textbf{Proof.} Considering of the proof of Theorem 5.1, we only need to show that $\mathcal{C}$ joins
$\left(\lambda_1/f_0, 0\right)$ to
$\left(0, +\infty\right)$. Corollary 5.1 implies that $\mathcal{C}$ is unbounded in the direction of $X^+$. Clearly, $f_\infty=+\infty$ implies that $f(s)\geq M^Ns^N$
for some positive constant $M$ and $s$ large enough.

To complete the proof, it suffices to show that the
unique blow up point of $\mathcal{C}$ is $\lambda=0$. Suppose on the contrary
that there exists $\widehat{\lambda}>0$ such that $\left(\widehat{\lambda},0\right)$ is a
blow up point of $\mathcal{C}$. Then there exists a sequence $\left\{\lambda_n, v_n\right\}$
such that $\underset{n\rightarrow +\infty}{\lim}\lambda_n=\widehat{\lambda}$ and
$\underset{n\rightarrow +\infty}{\lim}\left\Vert v_n\right\Vert=+\infty$. Let $w_n =v_n/\left\Vert v_n\right\Vert$. Clearly, one has that
\begin{equation}
w_n=\lambda_n\int_r^1\left(\int_0^sN\tau^{N-1}a(\tau)\left(\frac{f\left(v_n\right)}
{v_n^N}\frac{v_n^N}{\left\Vert v_n\right\Vert^N}\right)\,d\tau\right)^{1/N}\,ds.\nonumber
\end{equation}
Take $M=1/\left(\widehat{\lambda}\xi\right)+1$, where $\xi=\frac{1}{16}\left(\int_0^{1/4}N\tau^{N-1}a(\tau)\,d\tau\right)^{1/N}>0$. For $r\in\left[1/4,3/4\right]$, by virtue of Lemma 2.2
of [\ref{HW}], we have that
\begin{eqnarray}\label{C3}
\left\vert w_n\right\vert&\geq& M\lambda_n\int_r^1
\left(\int_0^sN\tau^{N-1}a(\tau)\left\vert w_n\right\vert^N\,d\tau\right)^{1/N}\,ds\nonumber\\
&\geq& M\left\Vert w_n\right\Vert\lambda_n\int_r^1
\left(\int_0^sN\tau^{N-1}a(\tau)(1-\tau)^N\,d\tau\right)^{1/N}\,ds\nonumber\\
&\geq& M\left\Vert w_n\right\Vert\lambda_n(1-r)
\left(\int_0^rN\tau^{N-1}a(\tau)(1-\tau)^N\,d\tau\right)^{1/N}\nonumber\\
&\geq& M\left\Vert w_n\right\Vert\lambda_n(1-r)^{2}
\left(\int_0^rN\tau^{N-1}a(\tau)\,d\tau\right)^{1/N}\nonumber\\
&\geq& M\left\Vert w_n\right\Vert\lambda_n \xi.
\end{eqnarray}
It is obvious that (\ref{C3}) follows $M\lambda_n\leq 1/\xi$. Thus, we
get that $M\leq 1/\left(\widehat{\lambda}\xi\right)$. While, this is impossible because
of $M=1/\left(\widehat{\lambda}\xi\right)+1$.\qed
\\ \\
\textbf{Remark 5.3.} Clearly, Theorem 5.3 and Corollary 5.1 imply that
if $f_0\in(0,+\infty)$ and $f_\infty=+\infty$ then there exists $\lambda_5>0$ such that
(\ref{MO}) has at least one strictly convex solution for all $\lambda\in\left(0,\lambda_5\right)$
and has no convex solution for all $\lambda\in\left(\lambda_5,+\infty\right)$.
\\ \\
\textbf{Theorem 5.4.} \emph{If $f_0=0$ and $f_\infty\in(0,+\infty)$,
then for any $\lambda\in\left(\lambda_1/f_\infty,+\infty\right)$, (\ref{MO})
has at least one solution $u$ such that it is negative, strictly convex in $(0,1)$.}
\\ \\
\textbf{Proof.} If $(\lambda,v)$ is any solution of (\ref{C1}) with
$\Vert v\Vert\not\equiv 0$, dividing (\ref{C1}) by $\Vert v\Vert^{2N}$ and
setting $w=v/\Vert v\Vert^2$ yield
\begin{equation}
\left\{
\begin{array}{l}
\left(\left(-w'(r)\right)^N\right)'=\lambda^NNr^{N-1}a(r)\left(\frac{f(v)}
{\Vert v\Vert^{2N}}\right),\,\, r\in(0,1),\\
w'(0)=w(1)=0.
\end{array}
\right.\label{C4}
\end{equation}
Define
\begin{eqnarray}
\widetilde{f}(w)=\left\{
\begin{array}{lll}
\Vert w\Vert^{2N}f\left(\frac{w}{\Vert w\Vert^2}\right),\,\,&\text{if}\, w\neq 0,\\
0,~~~~~~~~~~~~~~~~~~~\,\, &\text{if}\,\, w=0.
\end{array}
\right.\nonumber
\end{eqnarray}
\indent Clearly, (\ref{C4}) is equivalent to
\begin{equation}\label{C5}
\left\{
\begin{array}{l}
\left(\left(-w'(r)\right)^N\right)'=\lambda^Nr^{N-1}a(r)\widetilde{f}(w),\,\, r\in(0,1),\\
w'(0)=w(1)=0.
\end{array}
\right.
\end{equation}
It is obvious that $(\lambda,0)$ is always the solution of (\ref{C5}).
By the simple computation, we can show that $\widetilde{f}_0=f_\infty$
and $\widetilde{f}_\infty=f_0$.

Now applying Theorem 5.2 and the inversion $w\rightarrow w/\Vert w\Vert^2=v$,
we can achieve our conclusion.\qed\\
\\
\textbf{Remark 5.4.} Under the assumptions of Theorem 5.4, we note there exists
$\lambda_6>0$ such that
(\ref{MO}) has at least one strictly convex solution
for all $\lambda\in\left(\lambda_6,+\infty\right)$
and has no convex solution for all $\lambda\in\left(0,\lambda_6\right)$.\\
 \\
\textbf{Theorem 5.5.} \emph{If $f_0=0$ and $f_\infty=0$,
then there exists $\lambda_*>0$
such that for any $\lambda\in\left(\lambda_*,+\infty\right)$, (\ref{MO}) has at least two
solutions $u_1$ and $u_2$ such that they are negative, strictly convex in $(0,1)$.}
\\ \\
\textbf{Proof.} Define
\begin{eqnarray}
f^n(s)=\left\{
\begin{array}{lll}
\frac{1}{n}s^N,\,\,\quad\quad\quad\quad\quad\quad\quad\,
\quad\quad\quad\,\,\quad & s\in\left[0,\frac{1}{n}\right],\\
\left(f\left(\frac{2}{n}\right)-\frac{1}{n^{N+1}}\right)ns+\frac{2}{n^{N+1}}
-f\left(\frac{2}{n}\right),\,\, & s\in\left(\frac{1}{n},\frac{2}{n}\right),\\
f(s),\quad\quad\quad\quad\quad\quad\quad\quad\quad\quad\quad\quad\,\,& s\in\left[\frac{2}{n},+\infty\right).
\end{array}
\right.\nonumber
\end{eqnarray}
Now, consider the following problem
\begin{equation}
\left\{
\begin{array}{l}
\left(\left(u'(r)\right)^N\right)'=\lambda^NNr^{N-1}a(r) f^n(-u(r))\,\, \text{in}\,\, 0<r<1,\\
u'(0)=u(1)=0.
\end{array}
\right.\nonumber
\end{equation}
Clearly, we can see that $\lim_{n\rightarrow+\infty}f^n(s)=f(s)$, $f_0^n=1/n$ and $f_\infty^n=f_\infty=0$.
Theorem 5.2 implies
that there exists a sequence unbounded continua $\mathcal{C}_n$ emanating from $\left(n\lambda_1, 0\right)$
and joining to $\left(+\infty, +\infty\right)$.

Take $z^*=\left(+\infty, 0\right)$. Theorem 2.2 implies that there exists an unbounded component $\mathcal{C}$ of $\liminf_{n\rightarrow +\infty}\mathcal{C}_n$ such that $z^*\in \mathcal{C}$.

We claim that $\mathcal{C}\cap([0,+\infty)\times\{0\})=\emptyset$.
Otherwise, there exists a sequence $\left\{\mu_n, v_n\right\}$ such that
$\underset{n\rightarrow +\infty}{\lim}\mu_n=\mu_{M}$ and
$\underset{n\rightarrow +\infty}{\lim}\left\Vert v_n\right\Vert=0$ as
$n\rightarrow+\infty$. Let $w_n =v_n/\left\Vert v_n\right\Vert$ and $w_n$
should be the solutions of the following problem
\begin{equation}
w=\mu_n\int_r^1\left(\int_0^sN\tau^{N-1}a(\tau)\left(\frac{f\left(v_n\right)}
{\left\Vert v_n\right\Vert^N}\right)\,d\tau\right)^{1/N}\,ds.\nonumber
\end{equation}
Similar to (\ref{egn0}), we can show that
\begin{equation}
\lim_{n\rightarrow+\infty}\frac{ f\left(v_n(r)\right)}{\left\Vert v_n\right\Vert^{N}}=0.\nonumber
\end{equation}
By the compactness of $T_f$, we obtain that for some convenient subsequence
$w_n\rightarrow w_0$ as $n\rightarrow+\infty$. Letting $n\rightarrow+\infty$,
we obtain that
$w_0\equiv0$. This contradicts $\left\Vert w_0\right\Vert=1$.

Next, we show that $\mathcal{C}$ joins $z^*$ to $\left(+\infty, +\infty\right)$.
We first show that $\mathcal{C}$ is unbounded in the direction of $X^+$.
Suppose on the contrary that there exists $M>0$ such that $\Vert v\Vert\leq M$ for any $(\lambda,v)\in \mathcal{C}$.
Let $w =v/\left\Vert v\right\Vert$. Obviously, one has that
\begin{eqnarray}
1\geq w&=&\lambda\int_r^1\left(\int_0^sN\tau^{N-1}a(\tau)\left(\frac{f\left(v\right)}
{\left\Vert v\right\Vert^N}\right)\,d\tau\right)^{1/N}\,ds\nonumber \\
&=&\lambda\int_r^1\left(\int_0^sN\tau^{N-1}a(\tau)\left(\frac{f\left(v\right)}
{v^N}\right)\left(\frac{v^N}
{\left\Vert v\right\Vert^N}\right)\,d\tau\right)^{1/N}\,ds. \nonumber
\end{eqnarray}
For any $\alpha\in\left(0,1/2\right)$, Lemma 2.2 of [\ref{HW}] shows that
\begin{eqnarray}
\min_{r\in[\alpha,1-\alpha]} v(r)\geq \alpha \Vert v\Vert. \nonumber
\end{eqnarray}
The signum condition (f1) and $\Vert v\Vert\leq M$ imply that there exists a positive
constant $\rho$ such that $f(v)/v^N\geq \min_{s\in[\alpha M,M]}f(s)/s^N\geq\rho$ for
any $r\in[\alpha,1-\alpha]$.
So for any $r\in[\alpha,1-\alpha]$, we have
\begin{eqnarray}
1\geq w(r)=\lambda \alpha \rho^{1/N}N^{1/N}\xi,\nonumber
\end{eqnarray}
where $\xi=\min_{r\in[\alpha,1-\alpha]}\int_r^1\left(\int_0^s\tau^{N-1}a(\tau)\,d\tau\right)^{\frac{1}{N}}\,ds>0$.
Hence, we have
\begin{eqnarray}
\lambda\leq \frac{1}{\alpha \rho^{1/N}N^{1/N}\xi}.\nonumber
\end{eqnarray}
This is a contradiction. By the same proof as in Theorem 5.2, we can show that the unique  blow up point of $\mathcal{C}$ is $\left(+\infty,0\right)$.
So $\mathcal{C}$ joins $(+\infty,0)$ to $\left(+\infty, +\infty\right)$.
\qed\\
\\
\textbf{Remark 5.5.} From Theorem 5.5 and Corollary 5.2, we can also see that there exists $\lambda_7>0$ such that
(\ref{MO}) has at least one strictly convex solution for all
$\lambda\in\left[\lambda_7,\lambda_*\right]$ and has no convex solution for all $\lambda\in\left(0,\lambda_7\right)$.
\\ \\
\textbf{Theorem 5.6.} \emph{If $f_0=0$ and $f_\infty=\infty$,
then for any $\lambda\in\left(0,+\infty\right)$, (\ref{MO}) has at least one solution $u$ such
that it is negative, strictly convex in $(0,1)$.}
\\ \\
\textbf{Proof.} Using an argument similar to that of Theorem 5.5, in view of Theorem 5.3, we can easily get the results of this theorem.\qed
\\ \\
\textbf{Theorem 5.7.} \emph{If $f_0=\infty$ and $f_\infty=0$,
then for any $\lambda\in\left(0,+\infty\right)$, (\ref{MO}) has at least one solution $u$ such
that it is negative, strictly convex in $(0,1)$.}
\\ \\
\textbf{Proof.} By an argument similar to that of Theorem 5.4 and the conclusions of Theorem 5.6, we can prove it.\qed
\\ \\
\textbf{Theorem 5.8.} \emph{If $f_0=\infty$ and $f_\infty\in(0,+\infty)$,
then for any $\lambda\in\left(0,\lambda_1/f_\infty\right)$, (\ref{MO}) has
at least one solution $u$ such that it is negative, strictly convex in $(0,1)$.}
\\ \\
\textbf{Proof.} By an argument similar to that of Theorem 5.4 and the conclusion of Theorem 5.3, we can obtain it.\qed
\\ \\
\textbf{Remark 5.6.} Similarly to Remark 5.3, there exists $\lambda_8>0$ such that
(\ref{MO}) has at least one strictly convex solution for all $\lambda\in\left(0,\lambda_8\right)$
and has no convex solution for all $\lambda\in\left(\lambda_8,+\infty\right)$.
\\ \\
\textbf{Theorem 5.9.} \emph{If $f_0=\infty$ and $f_\infty=\infty$,
then there exists $\lambda^*>0$
such that for any $\lambda\in\left(0,\lambda^*\right)$, (\ref{MO}) has at least two solutions $u_1$
and $u_2$ such that they are negative, strictly convex in $(0,1)$.}\\
\\
\textbf{Proof.} Define
\begin{eqnarray}
f^n(s)=\left\{
\begin{array}{lll}
ns^N,\,\,\quad\quad\quad\quad
\quad\quad\quad\quad\,\,\quad & s\in\left[0,\frac{1}{n}\right],\\
\left(f\left(\frac{2}{n}\right)-\frac{1}{n^{N-1}}\right)ns+\frac{2}{n^{N-1}}-f\left(\frac{2}{n}\right),\,\,& s\in
\left(\frac{1}{n},\frac{2}{n}\right),\\
f(s),\quad\quad\quad\quad\quad\,
\quad\quad\quad\quad\quad\,& s\in\left[\frac{2}{n},+\infty\right).
\end{array}
\right.\nonumber
\end{eqnarray}
Clearly, we can see that $\lim_{n\rightarrow+\infty}f^n(s)=f(s)$, $f_0^n=n$ and $f_\infty^n=f_\infty=\infty$.
Theorem 5.3 implies
that there exists a sequence unbounded continua $\mathcal{C}_n$ emanating from $\left(\lambda_1/n, 0\right)$
and joining to $\left(0, +\infty\right)$.

Taking $z^*=\left(0, 0\right)$, clearly $z^*\in \liminf_{n\rightarrow +\infty}\mathcal{C}_n$. The compactness of $T_f$ implies that $\left(\cup_{n=1}^{+\infty} \mathcal{C}_n\right)\cap B_R$ is pre-compact. Theorem 2.1 implies that $\mathcal{C}=\limsup_{n\rightarrow +\infty}C_n$ is unbounded closed connected such that $z^*\in \mathcal{C}$ and $\left(0, +\infty\right)\in \mathcal{C}$.
By an argument similar to that of Theorem 5.3, we can show that $\mathcal{C}\cap((0,+\infty)\times\{0\})=\emptyset$. \qed
\\ \\
\textbf{Remark 5.7.} By Theorem 5.9 and Corollary 5.1, we can see that there exists $\lambda_9>0$ such that
(\ref{MO}) has at least one strictly convex solution for all $\lambda\in\left[\lambda^*,\lambda_9\right]$
and has no convex solution for all
$\lambda\in\left(\lambda_9,+\infty\right)$. \\
\\
\textbf{Remark 5.8.} Clearly, the conclusions of Theorem 1.1 of [\ref{W}] and Theorem 5.1 of [\ref{HW}]
are the corollaries of Theorem 5.1--4.9.\\
\\
\textbf{Remark 5.9.} Let $f(s)=e^s$ and $a(x)\equiv 1$. It can be easily verified that $f_0=\infty$ and $f_\infty=\infty$.
This fact with Remark 5.7 implies that there is no solution of problem (\ref{MO}) with $\lambda$ large enough,
and for sufficiently small $\lambda$ there are two strictly convex solutions. Set $\mu:=\lambda^{1/2}$.
Through a scaling, we can show that problem (\ref{MO}) is equivalent to
\begin{eqnarray}\label{MAh}
\left\{
\begin{array}{lll}
\det\left(D^2u\right)=e^{-u}\,\, &\text{in}\,\, B_\mu(0),\\
u=0~~~~~~~~~~~~~~\,&\text{on}\,\, \partial B_\mu(0),
\end{array}
\right.
\end{eqnarray}
where $B_\mu(0)$ denotes the set of $\{x\in \mathbb{R}^N: \vert x\vert\leq \mu\}$.
Hence there is no solution of problem (\ref{MAh}) with $\mu$ large enough,
and for sufficiently small $\mu$ there are two strictly convex solutions. Obviously,
this result improves the corresponding one of [\ref{ZW}, Theorem 4.1]. So Theorem 4.1 of [\ref{ZW}] is
our corollary of Theorem 5.9.\\
\\
\textbf{Remark 5.10.} Obviously, the results of Theorem 5.1--4.9 are also valid on
$B_R(0)$ for any $R>0$.

\section{Exact multiplicity of convex solutions}

\quad\, In this section, under some more strict assumptions of $f$, we shall
show that the unbounded continuum which are obtained in Section 4
may be smooth curves. We just show the case of $f_0\in(0,+\infty)$ and $f_\infty=0$. Other cases are similar.

Firstly, we study the local structure of the bifurcation branch $\mathcal{C}$
near $\left(\lambda_1,0\right)$, which is obtained in Theorem 4.1.
Let $\mathbb{E}=\mathbb{R}\times X^+$, $\Phi(\lambda,v):=v-\lambda T_g(v)$ and
\begin{equation}
\mathcal{S}:=\overline{\left\{(\lambda,v)\in \mathbb{E}: \Phi(\lambda,v)=0,
v\neq0\right\}}^{\mathbb{E}}.\nonumber
\end{equation}
In order to formulate and prove main results of this section, it is convenient
to introduce some notations.
Given any $\lambda\in \mathbb{R}$ and $0 < s < +\infty$, we consider an
open neighborhood of $\left(\lambda_1, 0\right)$ in $\mathbb{E}$ defined by
\begin{equation}
\mathbb{B}_s(\lambda_1,0):=\left\{(\lambda,v)\in \mathbb{E}: \Vert v\Vert
+\left\vert\lambda-\lambda_1\right\vert<s\right\}.\nonumber
\end{equation}
Let $X_0$ be a closed subspace of $X$ such that
\begin{equation}
X=\text{span}\left\{\psi_1\right\}\oplus X_0.\nonumber
\end{equation}
According to the Hahn-Banach theorem,
there exists a linear functional $l\in \left(X^+\right)^*$, here $\left(X^+\right)^*$ denotes the dual
space of $X^+$, such that
\begin{equation}
l\left(\psi_1\right)=1\,\, \text{and}\,\, X_0=\{v\in X^+: l(v)=0\}.\nonumber
\end{equation}
Finally, for any $0<\varepsilon<+\infty$ and $0 < \eta < 1$,
we define
\begin{equation}
K_{\varepsilon,\eta}^+:=\left\{(\lambda,v)\in \mathbb{E}: \left\vert\lambda-\lambda_1
\right\vert<\varepsilon, l(v)>\eta\Vert v\Vert\right\}.\nonumber
\end{equation}
\indent Applying an argument similar to that of [\ref{LG}, Lemma 6.4.1], we may obtain the following result, which localizes
the possible solutions of (\ref{Mg})
bifurcating from $\left(\lambda_1,0\right)$.\\ \\
\textbf{Lemma 6.1.} \emph{For every $\eta\in(0, 1)$ there exists a number
$\delta_0>0$ such that
for each $0<\delta<\delta_0$,
\begin{equation}
\left(\left(\mathcal{S}\setminus\left\{\left(\lambda_1,0\right)\right\}\right)\cap \mathbb{B}_\delta
\left(\lambda_1,0\right)\right)\subset K_{\varepsilon,\eta}^+.\nonumber
\end{equation}
Moreover, for each
\begin{equation}
(\lambda,v)\in\left(\mathcal{S}\setminus\left\{\left(\lambda_1,0\right)\right\}\right)\cap
\left(\mathbb{B}_\delta\left(\lambda_1,0\right)\right),\nonumber
\end{equation}
there are $s\in\mathbb{R}$ and unique $y\in X_0$ such that
\begin{equation}
v=s\psi_1+y\,\, \text{and} \,\, s>\eta\Vert v\Vert.\nonumber
\end{equation}
Furthermore, for these solutions $(\lambda,v)$,
\begin{equation}
\lambda=\lambda_1+o(1)\,\, \text{and}\,\, y=o(s)\nonumber
\end{equation}
as $s\rightarrow 0^+$.}\\
\\
\textbf{Remark 6.1.} From Lemma 6.1, we can see that $\mathcal{C}$ near
$\left(\lambda_1,0\right)$ is given by a curve
$(\lambda(s),v(s))=\left(\lambda_1+o(1),s\psi_1+o(s)\right)$ for $s$ near $0^+$.
\\

\indent The primary result in this section is the following theorem.
\\ \\
\textbf{Theorem 6.1.} \emph{Let $f\in C^1(\mathbb{R}\setminus\mathbb{R}^-,\mathbb{R}\setminus\mathbb{R}^-)$ satisfy the
assumptions of Theorem 5.2. Suppose $f'(s)<Nf(s)/s$ for any $s>0$.
Then for any $\lambda\in\left(\lambda_1/f_0,+\infty\right)$, (\ref{MO})
has exactly one solution $u_\lambda$ such that it is negative, strictly convex in $(0,1)$ and $u_\lambda$ decreasing with respect to $\lambda$.}
\\ \\
\textbf{Remark 6.2.} Clearly, the assumption $f'(s)<Nf(s)/s$ for $s>0$ is equivalent
to $f(s)/s^N$ is decreasing for $s>0$.
\\ \\
\indent We use the stability properties to prove Theorem 6.1.
Let
\begin{equation}
Y:=\left\{v\in C^2(0,1): v'(0)=v(1)=0\right\}.\nonumber
\end{equation}
For any $\phi\in Y$ and convex solution $u$ of (\ref{MO}),
by some simple computations, we can show
that the linearized equation of (\ref{MO}) about $u$ at the direction $\phi$ is
\begin{equation}\label{EM1}
\left\{
\begin{array}{l}
\left(-\phi'\left(-v'\right)^{N-1}\right)'-\lambda^Nr^{N-1}a(r)f'(v)\phi=\frac{\mu}{N} \phi\,\, \text{in}
\,\,(0,1),\\
\phi'(0)=\phi(1)=0,
\end{array}
\right.
\end{equation}
where $v=-u$. Hence, the linear stability of a solution $u$ of (\ref{MO})
can be determined by the linearized eigenvalue problem (\ref{EM1}).
A solution $u$ of (\ref{MO}) is stable if all eigenvalues of (\ref{EM1})
are positive, otherwise it is unstable.
We define the \emph{Morse index} $M(u)$ of a solution $u$ of (\ref{MO})
to be the number of negative eigenvalues of (\ref{EM1}).
A solution $u$ of (\ref{MO}) is degenerate if $0$ is an eigenvalue of
(\ref{EM1}), otherwise it is non-degenerate. \\

\indent The following lemma is our main stability result for the negative steady state solution.\\ \\
\textbf{Lemma 6.2.} \emph{Suppose that $f$ satisfies the conditions of Theorem 6.1.
Then any negative solution $u$ of (\ref{MO}) is stable, hence, non-degenerate and Morse index $M(u)=0$.}
\\ \\
\textbf{Proof.}
Let $u$ be a negative solution of (\ref{MO}), and let $\left(\mu_1, \varphi_1\right)$ be
the corresponding principal eigen-pairs of (\ref{EM1}) with $\varphi_1>0$ in $(0,1)$. We
notice that $v:=-u$ and $\varphi_1$ satisfy the
equations
\begin{equation}\label{EM2}
\left\{
\begin{array}{l}
\left(\left(-v'(r)\right)^N\right)'-\lambda^NNr^{N-1}a(r) f(v(r))=0\,\, \text{in}
\,\,(0,1),\\
v'(0)=v(1)=0
\end{array}
\right.
\end{equation}
and
\begin{equation}\label{EM3}
\left\{
\begin{array}{l}
\left(-\varphi_1'\left(-v'\right)^{N-1}\right)'-\lambda^Nr^{N-1}a(r)f'(v)\phi_1=\frac{\mu_1}{N} \varphi_1\,\, \text{in}
\,\,(0,1),\\
\phi_1'(0)=\phi_1(1)=0.
\end{array}
\right.
\end{equation}
Multiplying (\ref{EM3}) by $-v$ and (\ref{EM2}) by $-\varphi_1$, subtracting and integrating, we obtain
\begin{equation}
\mu_1\int_0^1 \varphi_1 v\,dr=N\int_0^1 \lambda^N r^{N-1}a(r)\varphi_1\left(Nf(v)-f'(v)v\right)\,dr.\nonumber
\end{equation}
Since $v> 0$ and $\varphi_1> 0$ in $(0,1)$, then $\mu_1 > 0$ and the negative steady state solution
$u$ must be stable.\qed\\
\\
\textbf{Proof of Theorem 6.1.} Define $F : \mathbb{R} \times X^+ \rightarrow X^+$ by
\begin{equation}
F(\lambda,v)=\left(\left(-v'(r)\right)^N\right)'-\lambda^NNr^{N-1}a(r) f(v(r)),\nonumber
\end{equation}
where $v=-u$. From Lemma 6.2, we know that any convex solution $v$ of (\ref{MO}) is stable.
Therefore, at any solution $\left(\lambda^*, v^*\right)$, we can apply Implicit Function Theorem to
$F(\lambda, v) = 0$, and all the solutions of $F(\lambda, v) = 0$ near $(\lambda^*, v^*)$ are on a curve
$(\lambda, v(\lambda))$ with $\left\vert \lambda-\lambda^*\right\vert\leq\varepsilon$ for some small $\varepsilon > 0$.
Furthermore, by virtue of Remark 6.1, the unbounded continuum $\mathcal{C}$ is a curve,
which has been obtained from Theorem 5.2.

Since $u_\lambda$ is differentiable with respect to
$\lambda$ (as a consequence of Implicit Function Theorem), letting $v_\lambda=-u_\lambda$, then
$\frac{dv_\lambda}{d\lambda}$ satisfies
\begin{equation}
\left(\left(\left(-\frac{d v_\lambda}{d \lambda}\right)'(r)\right)\left(-\left(v_\lambda\right)'(r)\right)^{N-1}\right)'
=\lambda^Nr^{N-1}a(r) f'\left(v_\lambda\right)\frac{d v_\lambda}{d \lambda}
+N\lambda^{N-1}r^{N-1}a(r)f\left(v_\lambda\right).\nonumber
\end{equation}
By an argument similar to that of Lemma 6.2, we can show that
\begin{equation}
\int_0^1\lambda^{N-1}r^{N-1}a(r)\left(\lambda\left(f'\left(v_\lambda\right)v_\lambda-Nf\left(v_\lambda\right)\right)\frac{d v_\lambda}{d \lambda}
+Nf\left(v_\lambda\right)v_\lambda\right)\,dr=0.\nonumber
\end{equation}
Assumptions of $f$ imply $\frac{d v_\lambda}{d \lambda}\geq0$. Therefore, we have $\frac{d u_\lambda}{d \lambda}\leq0$.\qed
\\ \\
\textbf{Remark 6.4.} From Theorem 6.1, we can get that (\ref{MO}) has no convex solution
for all $\lambda\in\left(0,\lambda_1/f_0\right]$ and for $\lambda>\lambda_1/f_0$ there exists a unique
negative solution $u_\lambda$ with $M\left(u_\lambda\right)=0$.
In this sense, we get the optical interval
for the parameter $\lambda$ which ensures the existence of single strictly convex solution for (\ref{MO}) under the assumptions of Theorem 6.1.
\\ \\
\textbf{Remark 6.4.} Note that the results of Theorem 6.1 have extended the corresponding
results to [\ref{L1}, Proposition 3] in the case of $\Omega=B$.\\
\\
\textbf{Remark 6.5.} Clearly, the results of Theorem 5.3 are better than the corresponding
results of [\ref{HW}, Theorem 3.1] if we assume $f\in C^1\left(\mathbb{R}\setminus \mathbb{R}^-,\mathbb{R}\setminus \mathbb{R}^-\right)$ in the Theorem 3.1 of [\ref{HW}].
Moreover, we do not need $f$ is increasing.
\\ \\
\textbf{Proof.} It is sufficient to show that the assumption (3.9) of [\ref{HW}] implies $f'(s)<Nf(s)/s$ for $s>0$.
Luckily, for any $s>0$ and $t\in(0,1)$, by the assumption (3.9) of [\ref{HW}], we have
\begin{eqnarray}
f'(s)&=&\lim_{t\rightarrow 1}\frac{f(s)-f(ts)}{(1-t)s}\leq\lim_{t\rightarrow 1}\frac{f(s)-\left[\left(1+\eta\right)t\right]^Nf(s)}{(1-t)s}\nonumber\\
&<&\lim_{t\rightarrow 1}\frac{f(s)-t^Nf(s)}{(1-t)s}=\lim_{t\rightarrow 1}\frac{\left(1+t+\cdots+t^{N-1}\right)f(s)}{s}=\frac{Nf(s)}{s},\nonumber
\end{eqnarray}
where $\eta>0$ comes from the assumption (3.9) of [\ref{HW}].\qed

\section{Convex solutions on general domain}

\quad\, In this section, we extend the results in Section 5 to the general domain $\Omega$ by domain comparison method.
\\

Through out this section, we assume that
\\

(f2) $f:[0,+\infty)\rightarrow[0,+\infty)$ is $C^2$ and decreasing;

(f3) $f(s)>0$ for $s>0$;

(A1) $a\in C^2\left(\overline{\Omega}\right)$.
\\

We use sub-supersolution method to construct a solution by iteration in an arbitrary domain. Note that 0 is always a sup-solution of problem (\ref{MAh1}). So we only need to find a sub-solution.

By an argument similar to that of [\ref{ZW}, Lemma 3.2], we
may obtain the following lemma.\\ \\
\textbf{Lemma 7.1.} \emph{If we have a strictly convex function $u_*\in C^3\left(\overline{\Omega}\right)$, such that
$\det\left(D^2u_*\right)\geq \lambda^Na(x)f\left(-u_*\right)$ in $\Omega$ and $u_*\leq 0$ on $\partial\Omega$, then problem (\ref{MAh1}) has a convex solution $u$ in $\Omega$.}
\\

As an immediate consequence, we obtain the following comparison.
\\ \\
\textbf{Lemma 7.2.} \emph{Given two bounded convex domains $\Omega_1$ and $\Omega_2$ such that
$\Omega_1\subset \Omega_2$. If we have a convex solution $u$ of problem (\ref{MAh1}) in $\Omega_2$,
then there exists a convex solution $v$ of problem (\ref{MAh1}) in $\Omega_1$, or equivalently if there is no convex solution of problem (\ref{MAh1}) in $\Omega_1$, then there is no convex solution of problem (\ref{MAh1}) in $\Omega_2$.}
\\

Our main results are the following two theorems.
\\ \\
\textbf{Theorem 7.1.} \emph{Assume that} (A1), (f2) \emph{and} (f3) \emph{hold}.\\

(a) \emph{If $f_0\in(0,+\infty)$ and $f_\infty\in(0,+\infty)$, then there exist $\lambda_2>0$ and $\lambda_3>0$ such that
(\ref{MAh1}) has at least one convex solution for all $\lambda\in\left(\lambda_2,\lambda_3\right)$.}

(b) \emph{If $f_0\in(0,+\infty)$ and $f_\infty=0$, then there exists $\lambda_4>0$ such that
(\ref{MAh1}) has at least one convex
solution for all $\lambda\in\left(\lambda_4,+\infty\right)$.}

(c) \emph{If $f_0\in(0,+\infty)$ and $f_\infty=+\infty$, then there exists $\lambda_5>0$ such that
(\ref{MAh1}) has at least one convex solution for all $\lambda\in\left(0,\lambda_5\right)$.}

(d) \emph{If $f_0=0$ and $f_\infty\in(0,+\infty)$, then there exists
$\lambda_6>0$ such that
(\ref{MAh1}) has at least one convex solution
for all $\lambda\in\left(\lambda_6,+\infty\right)$.}

(e) \emph{If $f_0=0$ and $f_\infty=0$, then there exist $\lambda_7>0$ and $\lambda_*>0$ such that
(\ref{MAh1}) has at least two convex solutions for all
$\lambda\in\left(\lambda_*,+\infty\right)$, one convex solution for all
$\lambda\in\left[\lambda_7,\lambda_*\right]$.}

(f) \emph{If $f_0=0$ (or $+\infty$) and $f_\infty=+\infty$ (or $0$), then for any $\lambda\in\left(0,+\infty\right)$, (\ref{MAh1}) has one convex solution.}

(g) \emph{If $f_0=+\infty$ and $f_\infty\in(0,\infty)$, then there exists $\lambda_8>0$ such that
(\ref{MAh1}) has at least one convex solution
for all $\lambda\in\left(0,\lambda_8\right)$.}

(h) \emph{If $f_0=+\infty$ and $f_\infty=+\infty$, then there exist $\lambda_9>0$ and $\lambda^*>0$ such that
(\ref{MAh1}) has at least two convex solutions for all $\lambda\in\left(0,\lambda^*\right)$, has at least one convex solution for all $\lambda\in\left[\lambda^*,\lambda_9\right]$.}
\\ \\
\textbf{Proof.} We only give the proof of (a) since the proofs of (b)--(h) can be given similarly.
It is obvious that there exists a positive constant $R_1$ such that $\Omega\subseteq B_{R_1}(0)$.
Theorem 5.1, Remark 5.1 and 5.10 that there exist $\lambda_2>0$ and $\lambda_3>0$ such that
problem (\ref{MAh1}) with $\Omega=B_{R_1}(0)$ has at least a strictly convex solution for all $\lambda\in\left(\lambda_2,\lambda_3\right)$.
Using Lemma 7.2, we have that problem (\ref{MAh1}) has at least a convex solution for all $\lambda\in\left(\lambda_2,\lambda_3\right)$.\qed
\\ \\
\textbf{Theorem 7.2.} \emph{Assume that} (A1), (f2) \emph{and} (f3) \emph{hold}.\\

(a) \emph{If $f_0\in(0,+\infty)$ and $f_\infty\in(0,+\infty)$, then there exist $\mu_2>0$ and $\mu_3>0$ such that
(\ref{MAh1}) has no convex solution for all
$\lambda\in\left(0,\mu_2\right)\cup\left(\mu_3,+\infty\right)$.}

(b) \emph{If $f_0\in(0,+\infty)$ and $f_\infty=0$, then there exists $\mu_4>0$ such that
(\ref{MAh1}) has no convex solution for all $\lambda\in\left(0,\mu_4\right)$.}

(c) \emph{If $f_0\in(0,+\infty)$ and $f_\infty=+\infty$, then there exists $\mu_5>0$ such that
(\ref{MAh1}) has no convex solution for all $\lambda\in\left(\mu_5,+\infty\right)$.}

(d) \emph{If $f_0=0$ and $f_\infty\in(0,+\infty)$, then there exists
$\mu_6>0$ such that
(\ref{MAh1}) has no convex solution for all $\lambda\in\left(0,\mu_6\right)$.}

(e) \emph{If $f_0=0$ and $f_\infty=0$, then there exists $\mu_7>0$ such that
(\ref{MAh1}) has no convex solution for all $\lambda\in\left(0,\mu_7\right)$.}

(f) \emph{If $f_0=+\infty$ and $f_\infty\in(0,\infty)$, then there exists $\mu_8>0$ such that
(\ref{MAh1}) has no convex solution for all $\lambda\in\left(\mu_8,+\infty\right)$.}

(g) \emph{If $f_0=+\infty$ and $f_\infty=+\infty$, then there exists $\mu_9>0$ such that
(\ref{MAh1}) has no convex solution for all
$\lambda\in\left(\mu_9,+\infty\right)$.}
\\ \\
\textbf{Proof.} We also only give the proof of (a) since the proofs of (b)--(g) can be given similarly.
It is obvious that there exists a positive constant $R_2$ such that $B_{R_2}(0)\subseteq \Omega$.
Theorem 5.1, Remark 5.1 and 5.10 imply that there exist $\mu_2>0$ and $\mu_3>0$ such that
problem (\ref{MAh1}) with $\Omega=B_{R_2}(0)$ has no convex solution for all
$\lambda\in\left(0,\mu_2\right)\cup\left(\mu_3,+\infty\right)$.
Using Lemma 7.2 again, we have that problem (\ref{MAh1}) has no convex solution for all
$\lambda\in\left(0,\mu_2\right)\cup\left(\mu_3,+\infty\right)$.\qed
\indent 


\begin{thebibliography}{99}

\bibitem{}\label{A} R.A. Admas, Sobolev spaces, New-York, Academic Press, 1975.

\bibitem{}\label{ACD} A. Ambrosetti, R.M. Calahorrano and F.R. Dobarro, Global branching for discontinuous problems, Comment. Math. Univ. Carolin. 31 (1990), 213--222.

\bibitem{}\label{AT} C.J. Amick and R.E.L. Turner, A global branch of steady vortex rings, J. Rein. Angew. Math.
384 (1988), 1--23.

\bibitem{}\label{ADT} D. Arcoya, J.I. Diaz and L. Tello, $S$-shaped bifurcation branch in a quasilinear multivalued model arising in climatoloty, J. Differential Equations 150 (1998), 215--225.

\bibitem{}\label{Bre} H. Brezis, Operateurs Maximaux Monotone et Semigroup de Contractions dans les Espase de Hilbert,
Math. Studies, vol. 5, North-Holland, Amsterdam, 1973.

\bibitem{}\label{CNS} L. Caffarelli, L. Nirenberg and J. Spruck, The Dirichlet problem for nonlinear second-order
elliptic equations, Part I. Monge-Amp\`{e}re equation, Comm. Pure Appl. Math. 37 (1984), 369--402.

\bibitem{}\label{CY1} S.Y. Cheng and S.T. Yau, On the regularity of the solution of the $n$-dimensional Minkowski problem, Comm.
Pure Appl. Math.29, (1976), 495--516.

\bibitem{}\label{CY} S.Y. Cheng and S.T. Yao, On the regularity of the Monge-Amp\`{e}re equation
$\det (\partial^2 u/\partial x_i\partial x_j) =F(x,u)$, Comm. Pure Appl. Math. 30 (1977), 41--68.

\bibitem{}\label{CY2} S.Y. Cheng and S.T. Yao, The real Monge-Amp\`{e}re equations and affine fiat structures,
(Proc. of the 1980 Beijing Symp. on Differential Geometry and Differential Equations),
Ed. S.S. Cheng and W.T. Wu, Science Press Beijing 1982, Gordon and Breach, New-York, 1982.

\bibitem{}\label{C} F.H. Clarke, Optimization and Nonsmooth Analysis, Wiley, New York, 1983.

\bibitem{}\label{Dai} G. Dai, Global branching for discontinuous problems involving the $p$-Laplacian, E. J. Differential Equations 44 (2013), 1--10.

\bibitem{}\label{DML} G. Dai, R. Ma and Y. Lu, Bifurcation from infinity and nodal solutions of quasilinear problems without the signum condition, J. Math. Anal. Appl. 397 (2013), 119--123.

\bibitem{}\label{DMW} G. Dai, R. Ma and H. Wang, Eigenvalues, bifurcation and one-sign solutions for the periodic $p$-Laplacian. Commun. Pure Appl. Anal. 12(6)  (2013), 2839--2872.

\bibitem{}\label{Del} Ph. Delano, Radially symmetric boundary value problems for real and complex
elliptic Monge-Amp\`{e}re equations, J. Differential Equations 58 (1985), 318--344.

\bibitem{}\label{DPEM} M. Del Pino, M. Elgueta and R. Man\'{a}sevich, A homotopic deformation along $p$ of
a Leray-Schauder degree result and existence for $(\vert u'\vert^{p-2}u')'+f(t,u)=0$, $u(0)=u(T)=0$, $p>1$, J. Differential Equations 80 (1989), 1--13.

\bibitem{}\label{De} K. Deimling, Nonlinear Functional Analysis, Springer-Verlag, New-York, 1987.

\bibitem{}\label{E} L.C. Evans, Partial Differential Equations, AMS, Rhode Island, 1998.

\bibitem{}\label{GT} D. Gilbarg and N.S. Trudinger, Elliptic Partial Differential Equations of Second Order, Springer-Verlag, Berlin, Heidelberg, 2001.

\bibitem{}\label{Guan3} P. Guan, $C^2$ a priori estimates for degenerate Monge-Amp\`{e}re equations, Duke Math. J. 86 (1997), 323--346.

\bibitem{}\label{Guan0} P. Guan, Topics in Geometric Fully Nonlinear Equations, unpublished notes, 2002.

\bibitem{}\label{Guan4} B. Guan and P. Guan, Convex hypersurfaces of prescribed curvatures, Ann. of Math. 156 (2002), 655--673.

\bibitem{}\label{Guan5} P. Guan and Y.Y. Li, The Weyl problem with nonnegative Gauss curvature, J. Differential Geom. 39 (1994), 331--342.

\bibitem{}\label{Guan6} P. Guan and C.S. Lin, On the equation $det\left(u_{ij} +\delta_{ij}u\right) = u^p f(x)$ on $S^n$, preprint, NCTS in Tsing-Hua
University, 2000.

\bibitem{}\label{Guan1} P. Guan and E. Sawyer, Regularity of subelliptic Monge-Amp\`{e}re equations in the plane, Trans. Amer. Math. Soc. 361 (2009), 4581--4591.

\bibitem{}\label{Guan2} P. Guan, N.S. Trudinger and X.J. Wang, On the Dirichlet problem for degenerate Monge-Amp\`{e}re equations, Acta Math. 182 (1999), 87--104.

\bibitem{}\label{HW} S. Hu and H. Wang, Convex solutions of BVP arising from Monge-Amp\`{e}re equations,
Discrete and Contin. Dyn. Syst. 16 (2006), 705--720.

\bibitem{}\label{IO} T. Idogawa and M. \^{O}tani, The first eigenvalues of some abstact elliptic
operator, Funkcialaj Ekvacioj, 38 (1995), 1--9.

\bibitem{}\label{ILL} B. Im, E. Lee and Y.H. Lee, A global bifurcation phenomena for second
order singular boundary value problems,
J. Math. Anal. Appl. 308 (2005), 61--78.

\bibitem{}\label{Jankowski} T. Jankowski, Positive solutions to Sturm-Liouville problems with non-local boundary conditions, Proc. Roy. Soc. Edinburgh Sect. A  144 (2014), 119--138.

\bibitem{}\label{K} N.V. Krylov, On degenerate nonlinear elliptic equations, Mat. Sbornik, 120 (1983), 311--330.

\bibitem{}\label{K1} N.D. Kutev, Nontrivial solutions for the equations of Monge-Amp\`{e}re type, J.
Math. Anal. Appl. 132 (1988), 424--433.

\bibitem{}\label{LiSun} H. Li and J. Sun, Positive solutions of sublinear Sturm-Liouville problems with changing
sign nonlinearity, Comput. Math. Appl., 58 (2009), 1808--1815.

\bibitem{}\label{L} P.L. Lions, Sur les \'{e}quations de Monge-Amp\`{e}re, I,
Manuscripta Math. 41 (1983) 1-43; II, Arch. Rat. Mech. Anal. Announced in C.R. Acad. Sci. Paris, 293 (1981), 589-592.

\bibitem{}\label{L1} P.L. Lions, Two remarks on Monge-Amp\`{e}re equations, Ann. Mat.
Pura Appl. (4) 142 (1985), 263--275.

\bibitem{}\label{LG} J. L\'{o}pez-G\'{o}mez, Spectral theory and nonlinear functional
analysis, Chapman and Hall/CRC, Boca
Raton, 2001.

\bibitem{}\label{MaAn} R. Ma and Y. An, Global structure of positive solutions for nonlocal boundary value
problems involving integral conditions, Nonlinear Anal. 71 (2009), 4364--4376.

\bibitem{}\label{MGX} R. Ma, C. Gao and J. Xu, Existence of positive solutions for first order discrete periodic boundary value problems with delay, Nonlinear Anal. 74(12) (2011), 4186--4191.

\bibitem{}\label{Minkowski} H. Minkowski, Allgemeine Lehrs\"{a}tze \"{u}ber die konvexen Polyeder, Nachr.Ges. Wiss. Gottingen, (1897), 198--219.

\bibitem{}\label{Nirenberg} L. Nirenberg, The Weyl and Minkowski problems in differential geometry in the large, Comm. Pure Applied
Math. 6 (1953), 337--394.

\bibitem{}\label{P3} A.V. Pogorelov, Extrinsic geometry of convex surfaces, Transl. Math. Monographs, 35AMS, Providence,
R.I., 1973.

\bibitem{}\label{P0} A. V. Pogorelov, Regularity of a convex surface with given Gaussian curvature, Mat. Sb. 31 (1952), 88--103.

\bibitem{}\label{P} A.V. Pogorelov, On the regularity of generalized solutions of the
equation $\det (\partial^2 u/\partial x_i\partial x_j) =\phi(x_1,\ldots,x_n)>0$,
Soviet Math. Dokl. 12 (1971), 1436--1440.

\bibitem{}\label{P1} A.V. Pogorelov, The Diriehlet problem for the $n$-dimensional
analogue of the Monge-Amp\`{e}re equation, Soviet. Math. Dokl. 12 (1971), 1727--1731.

\bibitem{}\label{P2} A.V. Pogorelov, The Minkowski multidimensional problem, J. Wiley, New-York, 1978.

\bibitem{}\label{R} P.H. Rabinowitz, Some global results for nonlinear eigenvalue problems,
J. Funct. Anal. 7 (1971), 487--513.

\bibitem{}\label{S} J. Shi, Blow up points of solution curves for a semilinear problem,
Topo. Meth. Nonl. Anal. 15 (2000), 251--266.

\bibitem{}\label{T} K. Tso, On a real Monge-Amp\`{e}re functional, Invent. Math. 101 (1990), 425--448.

\bibitem{}\label{Whyburn} G.T. Whyburn, Topological Analysis, Princeton University Press, Princeton, 1958.

\bibitem{}\label{W} H. Wang, Convex solutions of boundary value problems,
J. Math. Anal. Appl. 318 (2006), 246--252.

\bibitem{}\label{XOC} X. Xu, D. O'Regan and Y. Chen, Structure of positive solution sets of semi-positone singular boundary value problems, Nonlinear Anal. 72 (2010), 3535--3550.

\bibitem{}\label{XSO} X. Xu, J. Sun and D. O'Regan, Global structure of positive solution sets of nonlinear operator equations, Monatsh. Math. 165 (2012), 271--303.

\bibitem{}\label{ZW} Z. Zhang and K. Wang, Existence and non-existence of solutions for a class
of Monge-Amp\`{e}re equations, J. Differential Equations 246 (2009), 2849--2875.
\end{thebibliography}
\end{document}